\documentclass{article}

\usepackage{amsthm,amsmath,amsfonts,amssymb}
\usepackage{graphicx}
\usepackage{hyperref}

\numberwithin{equation}{section}

\def \half{\frac{1}{2}}

\def \demb {{{ \bigcap\kern -12.2pt{_\big\downarrow}}}}
\def \uemb {{{ \bigcup\kern -5.8pt{^\big\uparrow}}}}

\def \Re{{I \kern -.3em R}}
\def \Rn {{{\bf I\kern -1.6pt{\bf R}}}^{\rm n}}

\def \half {{1\over 2}}

\newtheorem{theorem}{Theorem}[section]
\newtheorem{lemma}[theorem]{Lemma}
\newtheorem{corollary}[theorem]{Corollary}
\newtheorem{remark}[theorem]{Remark}
\newtheorem{definition}[theorem]{Definition}
\newtheorem{proposition}[theorem]{Proposition}

\begin{document}

\title{Continuity of the set equilibria of non-autonomous damped wave equations with terms concentrating on the boundary}

\author{Gleiciane da Silva Arag\~ao\thanks{Departamento de Ci\^encias Exatas e da Terra, Universidade Federal de S\~ao Paulo. Av. Concei\c{c}\~ao, 515, Centro, Cep 09920-000, Diadema-SP, Brazil, Phone: +55 (11) 4044-0500. Email: gleiciane.aragao@unifesp.br. Partially supported by CNPq 475146/2013-1, Brazil.} \thanks{Corresponding author.} \, and \,  Flank David Morais Bezerra\thanks{Departamento de Matem\'atica, Universidade Federal da Para\'iba, Cidade Universit\'aria, Campus I, Via Expressa Padre Z\'e-Castelo Branco III, Cep 58051-900, Jo\~ao Pessoa-PB, Brazil. Phone: +55 (83) 3216-7434. Email: flank@mat.ufpb.br.}}

\date{}

\maketitle \thispagestyle{empty} \vspace{-10pt}

\begin{abstract}
In this paper we are interested in the behavior of the solutions of non-autonomous damped wave equations when some reaction terms are concentrated in a neighborhood of the boundary and this neighborhood  shrinks to boundary as a parameter $\varepsilon$ goes to zero. We prove the continuity of the set equilibria of these equations. Moreover, if an equilibrium solution of the limit problem is hyperbolic, then we show that the perturbed equation has one and only one equilibrium solution nearby.
\end{abstract}

\vspace{0.2cm}


\vspace{0.2cm}



\section{Introduction}

\begin{figure}[!h]
  \centering
  \includegraphics[width=.30\columnwidth]{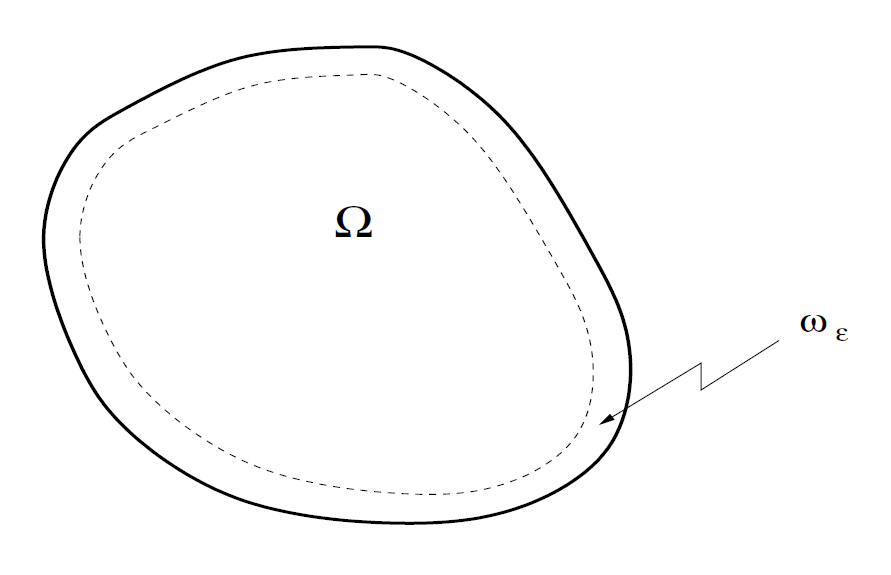} 
  \caption{The set $\omega_{\epsilon}$.}
  \label{figomega} 
\end{figure}

In this paper we prove the continuity of the set equilibria  of  non-autonomous damped wave equations when some reaction terms are concentrated in a neighborhood of the boundary and this neighborhood  shrinks to boundary as a parameter $\varepsilon$ goes to zero.  To better describe the problem, let $\Omega$ be an open bounded smooth set in  $\mathbb{R}^{3}$ with a smooth boundary $\Gamma=\partial \Omega$. We define the strip of width $\varepsilon$ and base $\partial \Omega$ as
$$
\omega_{\varepsilon}=\{x-\sigma\stackrel{\rightarrow}{n}(x): \ \mbox{$x\in \Gamma$  \ and  \ $\sigma \in [0,\varepsilon)$}\}, 
$$
for sufficiently small $\varepsilon$, say  $\varepsilon \in (0, \varepsilon_{0}]$, where $\stackrel{\rightarrow}{n}(x)$ denotes the outward normal vector at $x\in \Gamma$. We note that the set $\omega_{\varepsilon}$ has Lebesgue measure $\left|\omega_{\varepsilon}\right|=O(\varepsilon)$ with $\left|\omega_{\varepsilon}\right|\leqslant k\left|\Gamma\right|\varepsilon$, for some $k> 0$ independent of $\varepsilon$, and that for small $\varepsilon$, the set $\omega_{\varepsilon}$ is a neighborhood of $\Gamma$ in $\bar{\Omega}$, that collapses to the boundary when the parameter $\varepsilon$ goes to zero,  see Figure \ref{figomega}. 

We are interested in the behavior, for small $\varepsilon$, of the solutions of the non-autonomous damped wave equation with concentrated terms given by
\begin{equation}
\label{PPrin_1}
\begin{cases}
u_{tt}^\varepsilon-\mathrm{div}(a(x)\nabla u^\varepsilon)+u^\varepsilon+\beta(t)u_t^\varepsilon=f(u^\varepsilon)+\dfrac{1}{\varepsilon}\chi_{\omega_{\varepsilon}}g(u^\varepsilon)& \mbox{in}\ \Omega\times(\tau,+\infty),\\
a(x)\dfrac{\partial u^\varepsilon}{\partial\vec{n}}=0 & \mbox{on}\ \Gamma\times(\tau,+\infty),\\
u^\varepsilon(\tau)=u_0\in H^1(\Omega),\ \ u_t^\varepsilon(\tau)=v_0\in L^2(\Omega),
\end{cases}
\end{equation}
where $a\in\mathcal{C}^1(\overline{\Omega})$ with 
\[
0< a_0\leqslant a(x)\leqslant a_1,\quad \forall x\in\overline{\Omega},
\]
for fixed constants $a_0,a_1 >0$, and $\chi_{\omega_{\varepsilon}}$ denotes the characteristic function of the set $\omega_{\varepsilon}$. We refer to $\frac{1}{\varepsilon}\chi_{\omega_{\varepsilon}}g(u^\varepsilon)$ as the concentrating reaction  in $\omega_{\varepsilon}$. We assume that $\beta:\mathbb{R}\to\mathbb{R}$ is bounded, globally Lipschitz, and
\[
0<\beta_0\leqslant\beta(t)\leqslant\beta_1, \quad \forall t\in \mathbb{R},  
\]
for fixed constants $\beta_0,\beta_1> 0$ (the assumption that $\beta$ is globally Lipschitz continuity can be weakened to uniform continuity on $\mathbb{R}$ and continuous differentiability).

In \cite{aragaobezerra1} we take $f,g:\mathbb{R}\to\mathbb{R}$ to be $\mathcal{C}^2$ and assume that it satisfies the growth estimates
\begin{equation}
\label{Dcrescimento}
|j'(s)|\leqslant c(1+|s|^{\rho_j}),\quad \forall  s \in\mathbb{R},
\end{equation}
and
\begin{equation}
\label{crescimento}
|j(s_1)-j(s_2)|\leqslant c|s_1-s_2|(1+|s_1|^{\rho_j}+|s_2|^{\rho_j}),\quad \forall  s_1,s_2\in\mathbb{R},
\end{equation}
with  $j=f$ or $j=g$ and expoents $\rho_f$ and $\rho_g$, respectively, such that $\rho_f\leqslant 2$ and $\rho_g\leqslant 1$. We note that the estimate (\ref{Dcrescimento}) implies (\ref{crescimento}).

Moreover, we assume that
\begin{equation}
\label{DissipCond}
\limsup_{|s|\to+\infty}\dfrac{j(s)}{s}\leqslant0,
\end{equation}
with  $j=f$ or $j=g$. We note that \eqref{DissipCond} is equivalent to saying that for any $\gamma >0$ there exists $c_{\gamma}>0$ such that
\[
sj(s)\leqslant \gamma s^2+c_{\gamma}, \quad \forall s\in \mathbb{R}.
\]

As in \eqref{PPrin_1} the nonlinear term $g(u^\varepsilon)$ is only effective on the region $\omega_{\varepsilon}$ which collapses to $\Gamma$ as $\varepsilon\to0 $, then it is reasonable to expect that the family of solutions $u^{\varepsilon}$ of \eqref{PPrin_1}  will converge to a solution of an equation of the same type with nonlinear boundary condition on $\Gamma$. Indeed, under assumptions above, in \cite{aragaobezerra1}  we prove that the ``limit problem'' for the non-autonomous singularly wave equation \eqref{PPrin_1} is given by
\begin{equation}
\label{PPrin_2}
\begin{cases}
u_{tt}-\mathrm{div}(a(x)\nabla u)+u+\beta(t)u_t=f(u) & \mbox{in } \Omega\times(\tau,+\infty),\\
a(x)\dfrac{\partial u}{\partial\vec{n}}=g(u) & \mbox{on } \Gamma\times(\tau,+\infty),\\
u(\tau)=u_0\in H^1(\Omega),\ \ u_t(\tau)=v_0\in L^2(\Omega).
\end{cases}
\end{equation}
We also prove the existence and regularity of the pullback attractors of the problems (\ref{PPrin_1}) and (\ref{PPrin_2}), and that the family of attractors is upper semicontinuous at $\varepsilon=0$. Moreover, we show that the attractors are bounded in $H^2(\Omega)\times H^{1}(\Omega)$, uniformly in $\varepsilon$. In particular, all solutions of (\ref{PPrin_1}) and (\ref{PPrin_2}) are bounded in $H^{2}(\Omega)$ with a bound independent of $\varepsilon$ and, using Sobolev imbedding, we obtain that these solutions are also bounded in $L^{\infty}(\Omega)$, uniformly in $\varepsilon$. This enables us to cut the nonlinearities $f$ and $g$ in such a way that it becomes bounded with bounded derivatives up to second order without changing the solutions of the equations. After these considerations, without loss of generality, in all this work we will be assuming the following hypothesis 

\vspace{0.2cm}

\noindent \textbf{(H)} $f,g: \mathbb{R}\rightarrow \mathbb{R}$ are $\mathcal{C}^{2}$-functions satisfying (\ref{DissipCond}) and
$$
\left|j(u)\right|+\left|j'(u)\right|+\left|j''(u)\right|\leqslant K,\quad \mbox{$\forall u\in \mathbb{R}$},
$$
for some constant $K>0$ and with $j=f$ or $j=g$.

\vspace{0.2cm}

In this work we continue the analysis made in \cite{aragaobezerra1} and a first step to study the lower semicontinuity of the family of the pullback attractors of the problems (\ref{PPrin_1}) and (\ref{PPrin_2}) at $\varepsilon=0$ is to study the elements simpler of attractors, that is, the solutions of equilibria of (\ref{PPrin_1}) and (\ref{PPrin_2}). More precisely, we will prove the continuity of the set equilibria of \eqref{PPrin_1} and \eqref{PPrin_2} at $\varepsilon=0$. Also, we will prove an “uniqueness result,” in the sense that for any hyperbolic equilibrium of the limiting problem \eqref{PPrin_2}, there exists one and only one equilibrium of \eqref{PPrin_1} in its neighborhood. In particular, if all the equilibria of \eqref{PPrin_2} are hyperbolic, then there exists only a finite number of them and for all $\varepsilon$ small enough, the problem \eqref{PPrin_1} has exactly the same number of equilibria and they are close to the equilibria of \eqref{PPrin_2}.

In many theoretical and applied problems, it is important to understand what happens when the solutions varies parameters in the model, and wave equations with variable coefficients arise naturally in mathematical modeling of inhomogeneous media (e.g. functionally graded materials or materials with damage induced inhomogeneity) in solid mechanics, electromagnetic, fluid flows trough porous media (e.g. modeling  traveling waves in a inhomogeneous gas, see Egorov and Shubin \cite{YS} and Suggs \cite{Su}), and other areas of physics and engineering. Semilinear wave equation arises in quantum mechanics, whereas variants of the form $u_{tt}-\mathrm{div}(a\nabla u)+F(u,u_t)=0$
appear in the study of vibrating systems with or without damping, and with or without forcing terms; see e.g. Araruna and Bezerra \cite{AF} and Arrieta, Carvalho and Hale \cite{A-C-H}.  

Now, the non-autonomous damped wave equation
$$
u_{tt}-\Delta u+\beta(t)u_t=f(u), \quad \mbox{in $\Omega$}, 
$$
with Dirichlet boundary condition $u=0$ on $\Gamma$ and initial conditions $u(\tau)=u_0$ and $u_t(\tau)=v_0$ was studied in Carvalho, Langa and Robinson \cite[Chapter 15]{CLR}. Under the same assumptions as above, the authors used the theory of pullback asymptotic compactness to show that this equation had a pullback attractor and studied the gradient-like structure of the pullback attractor, see also Carabalho et al. \cite{CCLR_0}. 

On the other hand, concentrated terms problems on the strip $\omega_{\varepsilon}\subset \overline{\Omega}$ was initially studied in Arrieta,  Jim\'enez-Casas and Rodr\'iguez-Bernal  \cite{arrieta}, where linear elliptic equations with terms concentrated were considered and convergence results of the solutions were proved. Later, the asymptotic behavior of the attractors of a parabolic problem was analyzed in Jim\'enez-Casas and Rodr\'iguez-Bernal  \cite{anibal,AnibalAngela}, where the upper semicontinuity of attractors at $\varepsilon=0$ was proved. 
In these works the domain $\Omega$ is $C^2$ in $\mathbb{R}^n$, with $n\geqslant 2$. In \cite{gam1} some results of \cite{arrieta} were adapted to a nonlinear elliptic problem posed on an open square $\Omega$ in $\mathbb{R}^2$, considering $\omega_\epsilon \subset \Omega $ and with  highly oscillatory behavior in the boundary inside of $\Omega$. Later,  \cite{gam2} it proved the continuity of attractors for a nonlinear parabolic problem posed on a $C^2$ domain $\Omega$ in $\mathbb{R}^2$,  when some terms are concentrated in a neighborhood of the boundary and the inner boundary of this neighborhood presents a highly oscillatory behavior. The work \cite{gam2} was the first to consider the lower semicontinuity of attractors for problems with terms concentrating on the boundary. 

Note that in this works wave equations were not considered with concentration technique. However, in Arag\~ao and Bezerra \cite{aragaobezerra1} we considered a class of damped wave equations have Neumann boundary conditions and terms concentrating on the boundary and was proved the existence, regularity and upper semicontinuity of pullback attractors.  To our best knowledge, does not exist results  in the literature  on continuity of the set equilibria for these non-autonomous wave equations with terms concentrating on the boundary, and it is natural  because the hyperbolic structure of the equations here brings a further difficulty.

This paper is organized as follows. In Section \ref{SecAbstract}, we will give the notations that it will be used in this paper and we will see that the solutions of the nonlinear elliptic equations associated to the problems (\ref{PPrin_1}) and (\ref{PPrin_2}) will be obtained as fixed points of appropriate nonlinear maps defined in the space $H^{1}(\Omega)$. In Section \ref{SecTechnical}, we will prove several important technical results that will be needed in the proof of continuity of the set equilibria. In Section \ref{SecUpper}, we will show the upper semicontinuity of the set of equilibria of (\ref{PPrin_1}) and (\ref{PPrin_2}) at $\varepsilon=0$. Finally, in Section \ref{SecLower}, we will prove the lower semicontinuity of the set of equilibria of (\ref{PPrin_1}) and (\ref{PPrin_2}) at $\varepsilon=0$ and so we will obtain the continuity of this set of equilibria, for this, we will also need to assume that the equilibrium points of (\ref{PPrin_2}) are stable under perturbation. Moreover, if an equilibrium solution of the limit problem is hyperbolic, then we will show that the perturbed equation has one and only one equilibrium solution nearby.


\section{Abstract setting and solutions as fixed points}
\label{SecAbstract}

In this section we will present the functional framework that we will use to study the problems, we will define the abstract problems associated to (\ref{PPrin_1}) and (\ref{PPrin_2}). To better explain the results in the paper, initially, we will define the abstract problems associated to (\ref{PPrin_1}) and (\ref{PPrin_2}). Let us consider the  Hilbert space $Y:=L^2(\Omega)$, $Y^{\frac{1}{2}}:=H^1(\Omega)$ equipped with the inner product
\[
\langle u,v \rangle_{Y^{\frac{1}{2}}}=\int_\Omega a(x)\nabla u\nabla v dx+\int_\Omega uv dx,
\]
and we consider the unbounded linear operator $\Lambda:D(\Lambda)\subset Y\to Y$, defined by 
$$
\Lambda u=-\mbox{div}(a(x)\nabla u)+u, \quad u\in D(\Lambda),
$$
with
$$
D(\Lambda):=\Big\{u\in H^2(\Omega) \;:\; a(x)\dfrac{\partial u}{\partial\vec{n}}=0\ \mbox{on}\ \Gamma\Big\}.
$$
Since this operator turns out to be sectorial in $Y$, associated to it there is a scale of Banach spaces (the fractional power spaces)  $Y^\alpha$, $\alpha \in \mathbb{R}$, denoting the domain of the fractional power operators associated with $\Lambda$, that is, $Y^\alpha:=D(\Lambda^\alpha)$. Let us consider $Y^\alpha$ endowed with the graph norm $\|x\|_{Y^\alpha}=\|\Lambda^\alpha x\|_Y$ for $\alpha\geqslant0$. The fractional power spaces are related to the Bessel Potentials spaces $H^{s}(\Omega)$, $s\in \mathbb{R}$,  and
it is well know that 
\[
Y^\alpha \hookrightarrow H^{2\alpha}(\Omega),\quad Y^{-\alpha}=(Y^\alpha)',\ \ \alpha\geqslant0,
\]
with $Y^{\frac{1}{2}}=H^1(\Omega)$, $Y^{-\frac{1}{2}}=(H^1(\Omega))'$, $Y=Y^0=L^2(\Omega)$ and $Y^1=D(\Lambda)$. Since the problem (\ref{PPrin_2}) has a nonlinear term on boundary, choosing $\half < s \leqslant 1$ and using the standard trace theory results that for any function $v \in H^{s}(\Omega)$, the trace of $v$ is well defined and lies in $L^{2}(\Gamma)$. Moreover, the scale of negative exponents $Y^{-\alpha}$, for $\alpha > 0$, is necessary to introduce the nonlinear term of (\ref{PPrin_2}) in the abstract equation, since we are using the operator  $\Lambda$ with homogeneous boundary conditions.  Considering the realizations of $\Lambda$ in this scale, the operator $\Lambda_{-\frac{1}{2}} \in \mathcal{L}(Y^{\frac{1}{2}},Y^{-\frac{1}{2}})$ is given by 
$$
\langle  \Lambda_{-\half}u, v\rangle = \int_\Omega a(x)\nabla u\nabla vdx+\int_\Omega uv dx, \quad \mbox{for $u,v \in Y^{\frac{1}{2}}$.}
$$
With some abuse of notation we will identify all different realizations of this operator and we will write them all as $\Lambda$.

Also, let us consider the following Hilbert space
\[
X=X^0=Y^{\frac{1}{2}}\times Y
\]
equipped with the inner product
\[
\Big\langle\begin{pmatrix}u_1\\ v_1\end{pmatrix},\begin{pmatrix}u_2\\ v_2\end{pmatrix}\Big\rangle_X =\langle u_1,u_2\rangle_{Y^{\frac{1}{2}}}+\langle v_1,v_2 \rangle_Y,
\]
where $\langle \cdot,\cdot \rangle_Y$ is the usual inner product in $L^{2}(\Omega)$.

We define the unbounded linear operator $A:D(A)\subset X\to X$ by
$$
A\begin{pmatrix}u\\ v\end{pmatrix}=\begin{pmatrix}0 & -I\\ \Lambda & 0\end{pmatrix} \begin{pmatrix}u\\ v\end{pmatrix}=\begin{pmatrix}-v\\ \Lambda u\end{pmatrix}, \quad \begin{pmatrix}u\\ v\end{pmatrix} \in D(A), 
$$
with
$$
D(A)=D(\Lambda)\times Y= Y^1\times Y.
$$
It is proved in \cite[Proposition 6.21]{CLR} that $A$ generates a strongly continuous semigroup in $X$.

For each $\varepsilon \in (0,\varepsilon_{0}]$, we write \eqref{PPrin_1} in an abstract form as
\begin{equation}
\label{AP1}
\begin{cases}
\displaystyle \dfrac{dw^{\varepsilon}}{dt}+A w^{\varepsilon}=F_\varepsilon(t,w^{\varepsilon}), \quad  t> \tau,\\
w^{\varepsilon}(\tau)=w_0,
\end{cases}
\end{equation}
with 
\[
w^{\varepsilon}=\begin{pmatrix}u^{\varepsilon}\\ u^{\varepsilon}_t\end{pmatrix}, \quad w_{0}=\begin{pmatrix}u_{0}\\ v_{0}\end{pmatrix} \in X
\]
and nonlinear map $F_\varepsilon(t,\cdot): X  \to  H^{1}(\Omega) \times H^{-s}(\Omega)$, with $\frac{1}{2}< s\leqslant 1$ and $t> \tau$, defined by
\begin{equation*}
F_\varepsilon (t,w)=\begin{pmatrix}0\\ -\beta_{\Omega}(t,v)+f_{\Omega}(u)+\dfrac{1}{\varepsilon}\chi_{\omega_{\varepsilon}}g_{\Omega}(u)\end{pmatrix}, \quad \mbox{for } w=\begin{pmatrix}u\\ v\end{pmatrix} \in X,
\end{equation*}
where $\beta_{\Omega}(t,\cdot) :L^2(\Omega)\to H^{-s}(\Omega)$ and $f_{\Omega},\frac{1}{\varepsilon}\chi_{\omega_{\varepsilon}}g_{\Omega} :H^1(\Omega)\to H^{-s}(\Omega)$ are the operators, respectively, given by 
\begin{equation}
\label{betainterior}
\langle   \beta_{\Omega}(t,v),\varphi \rangle = \int_{\Omega} \beta(t)v \varphi dx,  \quad \mbox{$\forall v\in L^{2}(\Omega)$  and $\forall \varphi \in H^{s}(\Omega)$},
\end{equation}
\begin{equation}
\label{finterior}
\langle  f_{\Omega}(u),\varphi \rangle =\int_{\Omega} f(u) \varphi dx,  \quad \mbox{$\forall u\in H^{1}(\Omega)$  and $\forall \varphi \in H^{s}(\Omega)$} 
\end{equation}
and
\begin{equation}
\label{ginterior}
\langle  \dfrac{1}{\varepsilon}\chi_{\omega_{\varepsilon}}g_{\Omega}(u),\varphi \rangle =\frac{1}{\varepsilon} \int_{\omega_{\varepsilon}} g(u) \varphi dx,
 \quad \mbox{$\forall u\in H^{1}(\Omega)$  and $\forall \varphi \in H^{s}(\Omega)$}.
\end{equation}

While the problem \eqref{PPrin_2} can be written in an abstract form as
\begin{equation} 
\label{AP1b}
\begin{cases}
\displaystyle \dfrac{dw}{dt}+A w=F_0(t,w),\quad t> \tau,\\
w(\tau)=w_0,
\end{cases}
\end{equation}
with 
\[
w=\begin{pmatrix}u\\ u_t\end{pmatrix}
\]
and nonlinear map $F_0(t,\cdot): X \to H^{1}(\Omega) \times H^{-s}(\Omega)$, with $\frac{1}{2}< s\leqslant 1$  and $t> \tau$,  defined by
\begin{equation*}
F_0(t,w) =\begin{pmatrix}0\\ -\beta_{\Omega}(t,v)+f_{\Omega}(u)+g_{\Gamma}(u)\end{pmatrix}, \quad \mbox{for } w=\begin{pmatrix}u\\ v\end{pmatrix} \in  X,
\end{equation*}
where $\beta_{\Omega}(t,\cdot)$ and $f_{\Omega}$ are defined in (\ref{betainterior}) and (\ref{finterior}), respectively, and $g_{\Gamma} :H^1(\Omega)\to H^{-s}(\Omega)$ is the operator given by 
\begin{equation}
\label{gfronteira}
\langle  g_{\Gamma}(u),\varphi \rangle =\int_{\Gamma} \gamma(g(u)) \gamma(\varphi) dS, \quad \mbox{$\forall u\in H^{1}(\Omega)$  and $\forall \varphi \in H^{s}(\Omega)$},   
\end{equation}
where $\gamma: H^{s}(\Omega)\to  L^{2}(\Gamma)$ is the trace operator.

In \cite{aragaobezerra1} we obtain that for each $\tau \in\mathbb{R}$ and $w_0\in X$, the problems (\ref{AP1}) and (\ref{AP1b}) are global well-posedness and, for each $\varepsilon \in [0,\varepsilon_0]$, we many define an evolution process $\{S^\varepsilon(t,\tau):t\geqslant \tau\}$ in $X$ by 
$$
S^\varepsilon(t,\tau)w_0=w^\varepsilon(t,\tau,w_0),\quad t\geqslant \tau,
$$
where $w^\varepsilon$ is the unique solution of (\ref{AP1}) and (\ref{AP1b}). We prove the existence of pullback attractor $\{\mathcal{A}^\varepsilon(t): t\in\mathbb{R}\}$ for  \eqref{AP1} and \eqref{AP1b} in $X=H^1(\Omega)\times L^2(\Omega)$. We also prove the regularity and upper semicontinuity of the pullback attractors at $\varepsilon=0$. Here, we will continue this analysis now showing the continuity of the set equilibria of (\ref{AP1}) and (\ref{AP1b}) at $\varepsilon=0$. 


The equilibrium solutions of \eqref{PPrin_1} and \eqref{PPrin_2} (or of \eqref{AP1} and \eqref{AP1b}) are solutions that independent of time. Thus, we will need to consider the following nonlinear elliptic problems
\begin{equation}
\label{eql1}
\begin{cases}
-\mathrm{div}(a(x)\nabla u_\varepsilon)+u_\varepsilon=f(u_\varepsilon)+\dfrac{1}{\varepsilon}\chi_{\omega_{\varepsilon}}g(u_\varepsilon)& \mbox{in}\  \Omega,\\
a(x)\dfrac{\partial u_\varepsilon}{\partial\vec{n}}=0 & \mbox{on } \Gamma,
\end{cases}
\end{equation}
and
\begin{equation}
\label{eql2}
\begin{cases}
-\mathrm{div}(a(x)\nabla u)+u=f(u) & \mbox{in } \Omega,\\
a(x)\dfrac{\partial u}{\partial\vec{n}}=g(u) & \mbox{on } \Gamma.
\end{cases}
\end{equation}

Initially, we will write the elliptic problems \eqref{eql1} and \eqref{eql2} in abstract forms. For each $\varepsilon \in (0,\varepsilon_{0}]$, we write \eqref{eql1} in an abstract form as
\begin{equation}
\label{E1}
\displaystyle \Lambda u_{\varepsilon}=\tilde{F}_\varepsilon(u_{\varepsilon}),
\end{equation}
with $\tilde{F}_\varepsilon= f_{\Omega}+\frac{1}{\varepsilon}\chi_{\omega_{\varepsilon}}g_{\Omega}$, where $f_{\Omega},\frac{1}{\varepsilon}\chi_{\omega_{\varepsilon}}g_{\Omega}:H^1(\Omega)\to H^{-s}(\Omega)$, $\frac{1}{2}< s< 1$, are defined in (\ref{finterior}) and \eqref{ginterior}, respectively. While the problem \eqref{eql2} can be written in an abstract form as
\begin{equation} 
\label{E2}
\displaystyle \Lambda u=\tilde{F}_0(u),
\end{equation}
with $\tilde{F}_0= f_{\Omega}+g_{\Gamma}$, where $f_{\Omega},g_{\Gamma}:H^1(\Omega)\to H^{-s}(\Omega)$, $\frac{1}{2}< s< 1$, are defined in (\ref{finterior}) and \eqref{gfronteira}, respectively.

In particular, for each $\varepsilon\in [0,\varepsilon_{0}]$, $u_\varepsilon$ is a solution of (\ref{eql1}) and (\ref{eql2}) if and only if $u_\varepsilon \in H^{1}(\Omega)$ satisfies 
$$
u_{\varepsilon}=\Lambda^{-1} \tilde{F}_{\varepsilon}(u_{\varepsilon}),
$$
that is, $u_{\varepsilon}$ is a fixed point of the nonlinear map $\Lambda^{-1} \tilde{F}_{\varepsilon}: H^{1}(\Omega)\to H^{1}(\Omega)$.

For each $\varepsilon\in [0,\varepsilon_{0}]$, we denote by $\mathcal{E}_{\varepsilon}$ the set of solutions of (\ref{eql1}) and (\ref{eql2}), that is, 
$$
\mathcal{E}_{\varepsilon}= \left\{ u_{\varepsilon}\in H^{1}(\Omega):\ \Lambda u_{\varepsilon}-\tilde{F}_{\varepsilon}(u_{\varepsilon})=0 \right\}.
$$
Moreover, we denote by $E_{\varepsilon}$ the set of equilibria of \eqref{PPrin_1} and \eqref{PPrin_2}, that is, 
$$
E_{\varepsilon}= \left\{e_{\varepsilon}=\begin{pmatrix}u_{\varepsilon}\\ 0\end{pmatrix}\in H^{1}(\Omega) \times L^{2}(\Omega):\ u_{\varepsilon} \in \mathcal{E}_{\varepsilon}\right\}.
$$

Next, we will see that the upper semicontinuity of the family of equilibria $\{E_{\varepsilon}\}_{\varepsilon\in[0,\varepsilon_0]}$ of \eqref{PPrin_1} and \eqref{PPrin_2} at $\varepsilon=0$ is an immediate consequence of upper semicontinuity of the family attractors $\{\mathcal{A}^\varepsilon(t):t\in\mathbb{R}\}_{\varepsilon\in[0,\varepsilon_0]}$ at $\varepsilon=0$. On the other hand, note that in order to obtain the lower semicontinuity of the family of equilibria $\{E_{\varepsilon}\}_{\varepsilon\in[0,\varepsilon_0]}$, it is enough to prove the lower semicontinuity of the family of solutions $\{\mathcal{E}_{\varepsilon}\}_{\varepsilon\in[0,\varepsilon_0]}$ of (\ref{eql1}) and (\ref{eql2}) at $\varepsilon=0$. Before, we will show some technical results.


\section{Some technical results}
\label{SecTechnical} 

Initially, we will show that each set $\mathcal{E}_{\varepsilon}$ is not empty and it is compact.

\begin{lemma}
\label{compacto}
Suppose that (H) holds. Then, for each $\varepsilon\in [0,\varepsilon_{0}]$ fixed, the set $\mathcal{E}_{\varepsilon}$ of the solutions of (\ref{eql1}) and (\ref{eql2}) is not empty. Moreover, $\mathcal{E}_{\varepsilon}$ is compact in $H^{1}(\Omega)$.
\end{lemma}
\noindent {\bf Proof. }
Initially, we note that the linear operator $\Lambda^{-1}:H^{-s}(\Omega)\to H^{2-s}(\Omega)$ is continuous and using the compact embedding of $H^{2-s}(\Omega)$ in $H^{1}(\Omega)$, with $2-s>1$, we get that $\Lambda^{-1}:H^{-s}(\Omega)\to H^{1}(\Omega)$ is compact. Now, for each $\varepsilon\in [0,\varepsilon_{0}]$ fixed, using the \cite[Lemma 2.1]{aragaobezerra1}, we have that if $B$ is bounded set in $H^{1}(\Omega)$ then $\tilde{F}_{\varepsilon}(B)$ is bounded set in $H^{-s}(\Omega)$. Hence, by compactness of $\Lambda^{-1}:H^{-s}(\Omega)\to H^{1}(\Omega)$, we get that $\Lambda^{-1}\tilde{F}_{\epsilon}: H^{1}(\Omega)\to H^{1}(\Omega)$ is compact.

Now, show that for each $\varepsilon\in [0,\varepsilon_{0}]$ fixed, the set $\mathcal{E}_{\varepsilon}$ of the solutions of (\ref{eql1}) and (\ref{eql2}) is not empty, that is, that the equations (\ref{eql1}) and (\ref{eql2}) have at least one solution in $H^{1}(\Omega)$, it is equivalent to show that the compact operator $\Lambda^{-1}\tilde{F}_{\varepsilon}: H^{1}(\Omega)\to H^{1}(\Omega)$ has at least one fixed point. 

From \cite[Lemma 2.1]{aragaobezerra1}, we have that there exists $k>0$ independent of $\varepsilon$ such that
$$
\|\tilde{F}_{\varepsilon}(u)\|_{H^{-s}(\Omega)}\leqslant k,\quad \mbox{$\forall u\in H^{1}(\Omega)$\quad and\quad $\varepsilon \in [0, \varepsilon_{0}]$}.
$$

We consider the closed ball $\bar{B}_{r}(0)$ in $H^{1}(\Omega)$, where $r=k\|\Lambda^{-1}\|_{\mathcal{L}(H^{-s}(\Omega),H^{1}(\Omega))}$. For each $u\in H^{1}(\Omega)$, we have
\begin{eqnarray}
\label{limitacao}
\|\Lambda^{-1}\tilde{F}_{\varepsilon}(u)\|_{H^{1}(\Omega)}\leqslant\|\Lambda^{-1}\|_{\mathcal{L}(H^{-s}(\Omega),H^{1}(\Omega))} \|\tilde{F}_{\varepsilon}(u)\|_{H^{-s}(\Omega)}\leqslant r.
\end{eqnarray}

Therefore, the compact operator $\Lambda^{-1}\tilde{F}_{\varepsilon}: H^{1}(\Omega)\to H^{1}(\Omega)$ takes $H^{1}(\Omega)$ in the ball $\bar{B}_{r}(0)$, in particular, $\Lambda^{-1}\tilde{F}_{\varepsilon}$ takes $\bar{B}_{r}(0)$ into itself. From Schauder Fixed Point Theorem, we obtain that $\Lambda^{-1}\tilde{F}_{\varepsilon}$ has at least one fixed point in $H^{1}(\Omega)$. Thus, the equations (\ref{eql1}) and (\ref{eql2}) have at least one solution in $H^{1}(\Omega)$.

Now, for each $\varepsilon\in [0,\varepsilon_{0}]$ fixed, we will prove that $\mathcal{E}_{\varepsilon}$ is compact in $H^{1}(\Omega)$. For each $\varepsilon\in [0,\varepsilon_{0}]$ fixed, let $\{u_{n}\}_{n\in \mathbb{N}}$ be a sequence in $\mathcal{E}_{\varepsilon}$, then $u_{n}=\Lambda^{-1}\tilde{F}_{\varepsilon}(u_{n})$, for all $n\in \mathbb{N}$. Similarly to (\ref{limitacao}), we get that $\{u_{n}\}_{n\in \mathbb{N}}$ is a bounded sequence in $H^{1}(\Omega)$. Thus, $\{\Lambda^{-1}\tilde{F}_{\varepsilon}(u_{n})\}_{n\in \mathbb{N}}$ has a convergent subsequence, that we will denote by $\{\Lambda^{-1}\tilde{F}_{\varepsilon}(u_{n_{k}})\}_{k\in \mathbb{N}}$, with limit $u\in H^{1}(\Omega)$, that is,
$$
\Lambda^{-1}\tilde{F}_{\varepsilon}(u_{n_{k}})\to u \quad \mbox{in}\quad H^{1}(\Omega),\quad \mbox{as $k\to \infty$}.
$$
Hence, $u_{n_{k}}\to u$ in $H^{1}(\Omega)$, as $k\to \infty$. 

By continuity of operator $\Lambda^{-1}\tilde{F}_{\varepsilon}:H^{1}(\Omega)\to H^{1}(\Omega)$, we get 
$$
\Lambda^{-1}\tilde{F}_{\varepsilon}(u_{n_{k}})\to \Lambda^{-1}\tilde{F}_{\varepsilon}(u)\quad \mbox{in}\quad H^{1}(\Omega),\quad \mbox{as $k\to \infty$}.
$$
By the uniqueness of the limit, $u=\Lambda^{-1}\tilde{F}_{\varepsilon}(u)$. Thus, $\Lambda u-\tilde{F}_{\varepsilon}(u)=0$ and $u\in \mathcal{E}_{\varepsilon}$. Therefore, $\mathcal{E}_{\varepsilon}$ is a compact set in $H^{1}(\Omega)$. \quad $\blacksquare$

\vspace{0.2cm}

For each $\varepsilon\in(0,\varepsilon_0]$, we define the maps $Df_{\Omega},\frac{1}{\varepsilon}\chi_{\omega_{\varepsilon}}Dg_{\Omega},Dg_{\Gamma}:H^1(\Omega)\to \mathcal{L}(H^1(\Omega),H^{-s}(\Omega))$, respectively, by 
\begin{equation}
\label{Dfinterior}
\langle  Df_{\Omega}(u) w,\varphi \rangle =\int_{\Omega} f'(u)w \varphi dx,  \quad \mbox{$\forall u,w\in H^{1}(\Omega)$  and $\forall \varphi \in H^{s}(\Omega)$}, 
\end{equation}
\begin{equation}
\label{Dginterior}
\langle  \dfrac{1}{\varepsilon}\chi_{\omega_{\varepsilon}}Dg_{\Omega}(u) w,\varphi \rangle =\frac{1}{\varepsilon} \int_{\omega_{\varepsilon}} g'(u)w \varphi dx,
 \quad \mbox{$\forall u,w\in H^{1}(\Omega)$  and $\forall \varphi \in H^{s}(\Omega)$}
\end{equation}
and
\begin{equation}
\label{Dgfronteira}
\langle  Dg_{\Gamma}(u)w,\varphi \rangle =\int_{\Gamma} \gamma(g'(u)w) \gamma(\varphi) dS, \quad \mbox{$\forall u,w\in H^{1}(\Omega)$  and $\forall \varphi \in H^{s}(\Omega)$},   
\end{equation}
where $\gamma: H^{s}(\Omega)\to  L^{2}(\Gamma)$ is the trace operator. 


Now, we prove a result of uniform boundedness and convergence of the Fr\'echet differential of the nonlinearity $\tilde{F}_{\varepsilon}$.

\begin{lemma} 
\label{resultsconvnonlinearity}
Suppose that (H) holds. We have:
\begin{enumerate}
\item[1.)] For each $\varepsilon \in [0, \varepsilon_{0}]$, the map $\tilde{F}_{\varepsilon}:H^{1}(\Omega)\to H^{-s}(\Omega)$ is Fr\'echet differentiable, uniformly in $\varepsilon$, and your Fr\'echet differentials are given by $D\tilde{F}_0= Df_{\Omega}+Dg_{\Gamma}$ and $D\tilde{F}_\varepsilon= Df_{\Omega}+\frac{1}{\varepsilon}\chi_{\omega_{\varepsilon}}Dg_{\Omega}$, for $\varepsilon\in (0,\varepsilon_{0}]$, where the maps $Df_{\Omega},\frac{1}{\varepsilon}\chi_{\omega_{\varepsilon}}Dg_{\Omega},Dg_{\Gamma}:H^1(\Omega)\to \mathcal{L}(H^1(\Omega),H^{-s}(\Omega))$ are given respectively by \eqref{Dfinterior}, \eqref{Dginterior} and \eqref{Dgfronteira}.

\item[2.)] For each $\varepsilon \in [0, \varepsilon_{0}]$, the map $D\tilde{F}_{\varepsilon}:H^{1}(\Omega)\to \mathcal{L}(H^{1}(\Omega),H^{-s}(\Omega))$ is globally Lipschitz, uniformly in $\varepsilon$.

\item[3.)] There exist $k>0$ independent of $\varepsilon$ such that
$$
\|D\tilde{F}_{\varepsilon}(u^{*})\|_{\mathcal{L}(H^{1}(\Omega),H^{-s}(\Omega))}\leqslant k, \quad \mbox{$\forall u^{*} \in H^{1}(\Omega)$ \quad and \quad $\varepsilon \in [0, \varepsilon_{0}]$}. 
$$

\item[4.)] For each $u^{*}\in H^{1}(\Omega)$, we have
$$
\|D\tilde{F}_{\varepsilon}(u^{*})-D\tilde{F}_{0}(u^{*})\|_{\mathcal{L}(H^{1}(\Omega),H^{-s}(\Omega))}\to 0, \quad \mbox{as $\varepsilon\to0$},
$$
and this limit is uniform for $u^{*}\in H^{1}(\Omega)$ such that $\|u^{*}\|_{H^{1}(\Omega)}\leqslant R$, for some $R>0$.

\item[5.)] If $u^{*}_{\varepsilon}\to u^{*}$ in $H^{1}(\Omega)$, as $\varepsilon\to 0$, then
$$
\|D\tilde{F}_{\varepsilon}(u^{*}_{\varepsilon})-D\tilde{F}_{0}(u^{*})\|_{\mathcal{L}(H^{1}(\Omega),H^{-s}(\Omega))}\to 0, \quad \mbox{as $\varepsilon\to0$}.
$$

\item[6.)] If $u^{*}_{\varepsilon}\to u^{*}$ in $H^{1}(\Omega)$, as $\varepsilon\to 0$, and $w_{\varepsilon}\to w$ in $H^{1}(\Omega)$, as $\varepsilon\to 0$, then
$$
\|D\tilde{F}_{\varepsilon}(u^{*}_{\varepsilon}) w_{\varepsilon}-D\tilde{F}_{0}(u^{*}) w\|_{H^{-s}(\Omega)}\to 0, \quad \mbox{as $\varepsilon\to0$}.
$$
\end{enumerate}
\end{lemma}
\noindent {\bf Proof. }
The items {\it 1.)} and {\it 2.)} are immediate consequence of \cite[Lemmas 2.2 and 2.3]{aragaobezerra1}.

\vspace{0.2cm}

\noindent {\it 3.)} For each $u^{*}\in H^{1}(\Omega)$ and $\varepsilon \in [0, \varepsilon_{0}]$, we have
$$
\|D\tilde{F}_{\varepsilon}(u^{*})\|_{\mathcal{L}(H^{1}(\Omega),H^{-s}(\Omega))}=\displaystyle \sup_{\begin{array}{c}
w \in H^{1}(\Omega)\\
\left\|w\right\|_{H^{1}(\Omega)}=1
\end{array}} 
\| D\tilde{F}_{\varepsilon}(u^{*}) w \|_{H^{-s}(\Omega)}.
$$
Note that,  for each $w\in H^{1}(\Omega)$,
$$
\begin{array}{lll}
\| D\tilde{F}_{\varepsilon}(u^{*}) w \|_{H^{-s}(\Omega)}=\left\|Df_{\Omega}(u^{*}) w+\dfrac{1}{\varepsilon}\chi_{\omega_{\varepsilon}}Dg_{\Omega}(u^{*}) w\right\|_{H^{-s}(\Omega)}, \quad \varepsilon \in (0,\varepsilon_0], \\
\\
\| D\tilde{F}_{0}(u^{*}) w \|_{H^{-s}(\Omega)}=\|Df_{\Omega}(u^{*}) w+Dg_{\Gamma}(u^{*}) w\|_{H^{-s}(\Omega)},
\end{array}
$$
where the maps $Df_{\Omega},\frac{1}{\varepsilon}\chi_{\omega_{\varepsilon}}Dg_{\Omega}$ and $Dg_{\Gamma}$ are given respectively by \eqref{Dfinterior}, \eqref{Dginterior} and \eqref{Dgfronteira}.

Using that $f'$ is bounded and H\"older's inequality, we get 
\[
\begin{split}
\left|\langle Df_{\Omega}(u^{*}) w,\varphi\rangle \right| & \leqslant \int_\Omega |f'(u^{*})w||\varphi|dx \leqslant K\int_\Omega |w||\varphi|dx\\
&\leqslant K \Big[\int_\Omega|w|^2dx\Big]^{\frac{1}{2}}  \Big[\int_\Omega|\varphi|^2dx\Big]^{\frac{1}{2}}\\
&\leqslant  k_1 \|w\|_{H^1(\Omega)} \|\varphi\|_{H^s(\Omega)},\quad \forall \varphi \in H^{s}(\Omega).
\end{split}
\]
Thus, 
\begin{equation}
\label{F1}
\|Df_{\Omega}(u^{*}) w\|_{H^{-s}(\Omega)}\leqslant k_1\|w\|_{H^1(\Omega)}, \quad  \forall w \in H^1(\Omega). 
\end{equation}

Using that $g'$ is bounded, H\"older's inequality  and \cite[Lemma 2.1]{arrieta},  we have 
\[
\begin{split}
\left|\langle \dfrac{1}{\varepsilon}\chi_{\omega_{\varepsilon}}Dg_{\Omega}(u^{*}) w,\varphi\rangle \right|&\leqslant  \dfrac{1}{\varepsilon}\int_{\omega_{\varepsilon}} |g'(u^{*})w||\varphi|dx  \leqslant  \frac{K}{\varepsilon}\int_{\omega_{\varepsilon}}|w||\varphi|dx\\
&\leqslant K \Big[\dfrac{1}{\varepsilon}\int_{\omega_{\varepsilon}}|w|^2dx\Big]^{\frac{1}{2}}  \Big[\dfrac{1}{\varepsilon}\int_{\omega_{\varepsilon}}|\varphi|^2dx\Big]^{\frac{1}{2}}\\
& \leqslant k_2 \|w\|_{H^1(\Omega)}  \|\varphi\|_{H^s(\Omega)},\quad \forall \varphi \in H^{s}(\Omega),
\end{split}
\]
where the positive constant $k_2$ is independent of $\varepsilon$. Thus,

\begin{equation}
\label{F2}
\Big \|\dfrac{1}{\varepsilon}\chi_{\omega_{\varepsilon}}Dg_{\Omega}(u^{*}) w \Big \|_{H^{-s}(\Omega)}\leqslant k_2 \|w\|_{H^1(\Omega)}, \quad \forall w \in H^1(\Omega).
\end{equation}

Now,  using that $g'$ is bounded, H\"older's inequality and trace theorems, we get
\[
\begin{split}
\left|\langle Dg_{\Gamma}(u^{*}) w,\varphi\rangle \right|&\leqslant  \int_\Gamma |\gamma(g' (u^{*})w)||\gamma(\varphi)|dS \leqslant  K\int_\Gamma |\gamma(w)||\gamma(\varphi)|dS\\
&\leqslant K \Big[\int_\Gamma|\gamma(w)|^2dS\Big]^{\frac{1}{2}} \Big[\int_\Gamma|\gamma(\varphi)|^2dS\Big]^{\frac{1}{2}}\\
&\leqslant k_3 \|w\|_{H^1(\Omega)} \|\varphi\|_{H^s(\Omega)},\quad \forall \varphi \in H^{s}(\Omega).
\end{split}
\]
Thus,
\begin{equation}
\label{F3}
\|Dg_{\Gamma}(u^{*}) w\|_{H^{-s}(\Omega)}\leqslant k_3 \|w\|_{H^1(\Omega)}, \quad \forall w \in H^1(\Omega).
\end{equation}

Therefore, the result follows from \eqref{F1}, \eqref{F2} and \eqref{F3}.

\vspace{0.2cm}

\noindent {\it 4.)} For each $u^{*}\in H^{1}(\Omega)$, notice that
$$
\|D\tilde{F}_{\varepsilon}(u^{*})-D\tilde{F}_{0}(u^{*})\|_{\mathcal{L}(H^{1}(\Omega),H^{-s}(\Omega))}= \left\|\frac{1}{\varepsilon}\chi_{\omega_{\varepsilon}}Dg_{\Omega}(u^{*})-Dg_{\Gamma}(u^{*})\right\|_{\mathcal{L}(H^{1}(\Omega),H^{-s}(\Omega))}.
$$

For each $u^{*},w\in H^{1}(\Omega)$, in \cite[Lemma 5.2]{AnibalAngela} has been proven that there exists $M(\varepsilon,R)>0$ with $M(\varepsilon,R)\to 0$ as $\varepsilon \to 0$ such that
$$
\begin{array}{lll}
\displaystyle \left| \langle \dfrac{1}{\varepsilon}\chi_{\omega_{\varepsilon}}Dg_{\Omega}(u^{*}) w- Dg_{\Gamma}(u^{*}) w,\varphi\rangle  \right|  &=& \displaystyle \left| \frac{1}{\varepsilon}\int_{\omega_{\varepsilon}}g'(u^{*})w\varphi dx -\int_{\Gamma}\gamma(g'(u^{*})w)\gamma(\varphi) dS \right| \\

&\leqslant &  \displaystyle M(\varepsilon,R) \left\| w \right\|_{H^{1}(\Omega)} \left\| \varphi \right\|_{H^{1}(\Omega)}, \forall w,\varphi \in H^{1}(\Omega).
\end{array}
$$
Thus, 
\begin{equation}
\label{F4}
\left\| \dfrac{1}{\varepsilon}\chi_{\omega_{\varepsilon}}Dg_{\Omega}(u^{*})- Dg_{\Gamma}(u^{*})\right\|_{\mathcal{L}(H^{1}(\Omega),H^{-1}(\Omega))}  \to 0, \quad \mbox{as $\varepsilon \to 0$},
\end{equation}
uniformly for $u^{*}\in H^{1}(\Omega)$ such that $\left\| u^{*} \right\|_{H^{1}(\Omega)}\leqslant R$.

Now, fix $\frac{1}{2} < s_0 < 1$ and $0<\theta<1$. Then for any $s$ such that $-1<-s<-s_0<-\frac{1}{2}$,  using \eqref{F2}, \eqref{F3} and interpolation we have
$$
\begin{array}{lll}
\left\| \dfrac{1}{\varepsilon}\chi_{\omega_{\varepsilon}}Dg_{\Omega}(u^{*}) w- Dg_{\Gamma}(u^{*}) w\right\|_{H^{-s}(\Omega)} \\

\leqslant \left\| \dfrac{1}{\varepsilon}\chi_{\omega_{\varepsilon}}Dg_{\Omega}(u^{*}) w- Dg_{\Gamma}(u^{*}) w\right\|^{\theta}_{H^{-s_0}(\Omega)}\left\| \dfrac{1}{\varepsilon}\chi_{\omega_{\varepsilon}}Dg_{\Omega}(u^{*}) w- Dg_{\Gamma}(u^{*}) w\right\|^{1-\theta}_{H^{-1}(\Omega)}\\

\leqslant  (k_2+k_3)^{\theta} \left\| \dfrac{1}{\varepsilon}\chi_{\omega_{\varepsilon}}Dg_{\Omega}(u^{*})- Dg_{\Gamma}(u^{*})\right\|^{1-\theta}_{\mathcal{L}(H^{1}(\Omega),H^{-1}(\Omega))}\|w\|_{H^{1}(\Omega)}, \forall w\in H^{1}(\Omega).
\end{array}
$$
Hence and using \eqref{F4}, we obtain
$$
\left\| \dfrac{1}{\varepsilon}\chi_{\omega_{\varepsilon}}Dg_{\Omega}(u^{*})- Dg_{\Gamma}(u^{*})\right\|_{\mathcal{L}(H^{1}(\Omega),H^{-s}(\Omega))}  \to 0, \quad \mbox{as $\varepsilon \to 0$},
$$
uniformly for $u^{*}\in H^{1}(\Omega)$ such that $\left\| u^{*} \right\|_{H^{1}(\Omega)}\leqslant R$.

\vspace{0.2cm}

\noindent {\it 5.)} From item {\it 2.)}, we have that there exists $L>0$ independent of $\varepsilon$ such that
$$
\begin{array}{lll}
\displaystyle \|D\tilde{F}_{\varepsilon}(u^{*}_{\varepsilon})-D\tilde{F}_{0}(u^{*})\|_{\mathcal{L}(H^{1}(\Omega),H^{-s}(\Omega))}\\

\leqslant  \displaystyle \|D\tilde{F}_{\varepsilon}(u^{*}_{\varepsilon})-D\tilde{F}_{\varepsilon}(u^{*})\|_{\mathcal{L}(H^{1}(\Omega),H^{-s}(\Omega))}+\|D\tilde{F}_{\varepsilon}(u^{*})-D\tilde{F}_{0}(u^{*})\|_{\mathcal{L}(H^{1}(\Omega),H^{-s}(\Omega))}\\

\leqslant \displaystyle  L\|u^{*}_{\varepsilon}-u^{*}\|_{H^{1}(\Omega)}+\|D\tilde{F}_{\varepsilon}(u^{*})-D\tilde{F}_{0}(u^{*})\|_{\mathcal{L}(H^{1}(\Omega),H^{-s}(\Omega))}\to 0, \quad \mbox{as $\varepsilon\to0$},
\end{array}
$$
where we also use the item {\it 4.)} and $u^{*}_{\varepsilon}\to u^{*}$ in $H^{1}(\Omega)$, as $\varepsilon\to 0$.

\vspace{0.2cm}

\noindent {\it 6.)} We take $u^{*}_{\varepsilon}\to u^{*}$ in $H^{1}(\Omega)$, as $\varepsilon\to 0$, and $w_{\varepsilon}\to w$ in $H^{1}(\Omega)$, as $\varepsilon\to 0$. Using the items {\it 3.)} and {\it 5.)}, we get
$$
\begin{array}{lll}
\displaystyle \|D\tilde{F}_{\varepsilon}(u^{*}_{\varepsilon}) w_{\varepsilon}-D\tilde{F}_{0}(u^{*}) w\|_{H^{-s}(\Omega)}\\

\leqslant \displaystyle  \|D\tilde{F}_{\varepsilon}(u^{*}_{\varepsilon}) w_{\varepsilon}-D\tilde{F}_{\varepsilon}(u^{*}_{\varepsilon}) w\|_{H^{-s}(\Omega)}+\|D\tilde{F}_{\varepsilon}(u^{*}_{\varepsilon}) w-D\tilde{F}_{0}(u^{*}) w\|_{H^{-s}(\Omega)}\\

\leqslant \displaystyle \|D\tilde{F}_{\varepsilon}(u^{*}_{\varepsilon})\|_{\mathcal{L}(H^{1}(\Omega),H^{-s}(\Omega))}\|w_{\varepsilon}-w\|_{H^{1}(\Omega)}\\

+\|D\tilde{F}_{\varepsilon}(u^{*}_{\varepsilon})-D\tilde{F}_{0}(u^{*})\|_{\mathcal{L}(H^{1}(\Omega),H^{-s}(\Omega))}\|w\|_{H^{1}(\Omega)}\\

\leqslant \displaystyle k\|w_{\varepsilon}-w\|_{H^{1}(\Omega)}+\|D\tilde{F}_{\varepsilon}(u^{*}_{\varepsilon})-D\tilde{F}_{0}(u^{*})\|_{\mathcal{L}(H^{1}(\Omega),H^{-s}(\Omega))}\|w\|_{H^{1}(\Omega)}\to0, 
\end{array} 
$$
as $\varepsilon\to0$. \quad $\blacksquare$


\section{Upper semicontinuity of the set equilibria}
\label{SecUpper} 

We will prove the upper semicontinuity of the family of equilibria $\{E_{\varepsilon}\}_{\varepsilon\in[0,\varepsilon_0]}$ of \eqref{PPrin_1} and \eqref{PPrin_2} at $\varepsilon=0$. For this, we will use the fact that $E_{\varepsilon} \subset \mathcal{A}^\varepsilon(t)$ and that the family of attractors $\{\mathcal{A}^\varepsilon(t):t\in\mathbb{R}\}_{\varepsilon\in[0,\varepsilon_0]}$ is upper semicontinuous at $\varepsilon=0$.

\begin{theorem}\label{Theo_upper_equilibrio}
Suppose that (H) holds. Then the family of equilibria $\{E_{\varepsilon}\}_{\varepsilon\in[0,\varepsilon_0]}$ of \eqref{PPrin_1} and \eqref{PPrin_2} is upper semicontinuous at $\varepsilon=0$.
\end{theorem}
\noindent {\bf Proof. }
We will prove that for any sequence of $\varepsilon \to 0$ and for any $e_{\varepsilon} \in E_{\varepsilon}$ we can extract a subsequence which converges to an element of $E_{0}$. From the upper semicontinuity of the attractors and using that $e_{\varepsilon} \in E_{\varepsilon} \subset \mathcal{A}^\varepsilon(t)$, we can extract a subsequence $e_{\varepsilon_k} \in E_{\varepsilon_k}$ with $\varepsilon_k \to 0$, as $k\to \infty$, and we obtain the existence of a $e_{0} \in \mathcal{A}^0(t)$ such that 
$$
\|e_{\varepsilon_k}-e_0\|_{X} \to 0, \quad \mbox{as $k\to \infty$.}
$$
We need to prove that $e_0 \in E_{0}$, that is, $S^{0}(t,\tau)e_0=e_0$, for any $t\geqslant \tau$.

We first observe that for any $t>\tau$,
$$
\|e_{\varepsilon_k}-S^{0}(t,\tau)e_0\|_{X}\leqslant \|e_{\varepsilon_k}-e_0\|_{X}+\|e_{0}-S^{0}(t,\tau)e_0\|_{X}\to \|e_{0}-S^{0}(t,\tau)e_0\|_{X}, 
$$
as $k\to \infty$. Moreover, for a fixed $\tau^{*}>\tau$ and for any $t\in (\tau, \tau^{*})$, we obtain
$$
\begin{array}{lll}
\|e_{\varepsilon_k}-S^{0}(t,\tau)e_0\|_{X}\\

=\|S^{\varepsilon_k}(t,\tau)e_{\varepsilon_k}-S^{0}(t,\tau)e_0\|_{X}\\

\leqslant \|S^{\varepsilon_k}(t,\tau)e_{\varepsilon_k}-S^{0}(t,\tau)e_{\varepsilon_k}\|_{X}+\|S^{0}(t,\tau)e_{\varepsilon_k}-S^{0}(t,\tau)e_{0}\|_{X} \to 0,\quad \mbox{as $k\to \infty$,}
\end{array}
$$
where we have used the continuity of processes given by \cite[Lemma 6.1]{aragaobezerra1} and that $\{S^{0}(t,\tau):t\geqslant \tau\} \subset \mathcal{L}(X)$. In particular, we have that for each $t\geqslant \tau$, $S^{0}(t,\tau)e_0=e_0$, which implies that $e_0\in E_0$. \quad $\blacksquare$


\section{Lower semicontinuity of the set equilibria}
\label{SecLower}

To obtain the lower semicontinuity of the family of equilibria $\{E_{\varepsilon}\}_{\varepsilon\in[0,\varepsilon_0]}$ of \eqref{PPrin_1} and \eqref{PPrin_2} at $\varepsilon=0$, we need to prove the lower semicontinuity of the family of solutions $\{\mathcal{E}_{\varepsilon}\}_{\varepsilon\in[0,\varepsilon_0]}$ of (\ref{eql1}) and (\ref{eql2}) at $\varepsilon=0$ and this proof requeres additional assumptions. We need to assume that the solutions of \eqref{eql2} are stable under perturbation. This stability under perturbation can be given by the hyperbolicity. 

\begin{definition}\label{def_hyperbolic}
For each $\varepsilon\in[0,\varepsilon_{0}]$, we say that the solution $u^{*}_{\varepsilon}$ of (\ref{eql1}) and (\ref{eql2}) is hyperbolic if the spectrum $\sigma(\Lambda-D\tilde{F}_{\varepsilon}(u^{*}_{\varepsilon}))$ is disjoint from the imaginary axis, that is, $\sigma(\Lambda-D\tilde{F}_{\varepsilon}(u^{*}_{\varepsilon}))\cap \imath \mathbb{R}=\emptyset$. 
\end{definition}

\begin{theorem} \label{theoremisolated}
Suppose that (H) holds. If $u^{*}_{0}$ is a solution of \eqref{eql2} which satisfies $0\notin \sigma(\Lambda-D\tilde{F}_{0}(u^{*}_{0}))$, then $u^{*}_{0}$ is isolated.
\end{theorem}
\noindent {\bf Proof. }
Since $0\notin \sigma(\Lambda-D\tilde{F}_{0}(u^{*}_{0}))$ then $0\in \rho(\Lambda-D\tilde{F}_{0}(u^{*}_{0}))$. Thus, there exists $C>0$ such that 
$$
\|(\Lambda-D\tilde{F}_{0}(u^{*}_{0}))^{-1}\|_{\mathcal{L}(H^{-s}(\Omega),H^{1}(\Omega))}\leqslant C.
$$

Now, we note that $u$ is a solution of (\ref{eql2}) if and only if
$$
\begin{array}{rcl}
0&=& \Lambda u-D\tilde{F}_{0}(u^{*}_{0}) u+D\tilde{F}_{0}(u^{*}_{0}) u-\tilde{F}_{0}(u)\\

&\Leftrightarrow & (\Lambda-D\tilde{F}_{0}(u^{*}_{0})) u=\tilde{F}_{0}(u)-D\tilde{F}_{0}(u^{*}_{0}) u\\

&\Leftrightarrow & u=(\Lambda-D\tilde{F}_{0}(u^{*}_{0}))^{-1}(\tilde{F}_{0}(u)-D\tilde{F}_{0}(u^{*}_{0}) u). 
\end{array}
$$
So, $u$ is a solution of (\ref{eql2}) if and only if $u$ is a fixed point of the map
$$
\begin{array}{ll}
\Phi: H^{1}(\Omega)\to H^{1}(\Omega) \\
\qquad \qquad u\mapsto \Phi(u)=(\Lambda-D\tilde{F}_{0}(u^{*}_{0}))^{-1}(\tilde{F}_{0}(u)-D\tilde{F}_{0}(u^{*}_{0})u).
\end{array}
$$

We will show that there exists $r>0$ such that $\Phi:\bar{B}_{r}(u^{*}_{0})\to \bar{B}_{r}(u^{*}_{0})$ is a contraction, where $\bar{B}_{r}(u^{*}_{0})$ is a closed ball in $H^{1}(\Omega)$ with center in $u^{*}_{0}$ and ray $r$. 

In fact, from item {\it 1.)} of Lemma~\ref{resultsconvnonlinearity} we have that there exists $\delta>0$ such that
$$
C\|\tilde{F}_{0}(u)-\tilde{F}_{0}(v)-D\tilde{F}_{0}(u^{*}_{0})(u-v)\|_{H^{-s}(\Omega)}
\leqslant \frac{1}{2}\left\|u-v\right\|_{H^{1}(\Omega)},\quad \mbox{for}\quad \left\|u-v\right\|_{H^{1}(\Omega)}\leqslant \delta.
$$

Taking $r=\frac{\delta}{2}$ and $u$, $v\in \bar{B}_{r}(u^{*}_{0})$, we have
$$ 
\begin{array}{rcl}
\displaystyle \|\Phi(u)-\Phi(v)\|_{H^{1}(\Omega)}&=&\displaystyle \|(\Lambda-D\tilde{F}_{0}(u^{*}_{0}))^{-1}[\tilde{F}_{0}(u)-\tilde{F}_{0}(v)-D\tilde{F}_{0}(u^{*}_{0})(u-v)]\|_{H^{1}(\Omega)}\\

&\leqslant& \displaystyle C\|\tilde{F}_{0}(u)-\tilde{F}_{0}(v)-D\tilde{F}_{0}(u^{*}_{0})(u-v)\|_{H^{-s}(\Omega)}\\

&\leqslant & \frac{1}{2} \left\|u-v\right\|_{H^{1}(\Omega)}.
\end{array}
$$
Thus, $\Phi$ is a contraction on the $\bar{B}_{r}(u^{*}_{0})$. Moreover, if $u\in \bar{B}_{r}(u^{*}_{0})$ then
$$
\|\Phi(u)-u^{*}_{0}\|_{H^{1}(\Omega)}= \|\Phi(u)-\Phi(u^{*}_{0})\|_{H^{1}(\Omega)}\leqslant \frac{1}{2} \left\|u-u^{*}_{0}\right\|_{H^{1}(\Omega)}\leqslant \frac{r}{2}<r.
$$
Hence, $\Phi\left(\bar{B}_{r}(u^{*}_{0})\right)\subset \bar{B}_{r}(u^{*}_{0})$.

Therefore, from Contraction Theorem, $\Phi$ has an unique fixed point in $\bar{B}_{r}(u^{*}_{0})$. Since $u^{*}_{0}$ is a fixed point of $\Phi$, then $u^{*}_{0}$ is the unique fixed point of $\Phi$ in $\bar{B}_{r}(u^{*}_{0})$. Thus, $u^{*}_{0}$ is isolated.  \quad $\blacksquare$

\begin{corollary}
\label{corollaryisolated}
Suppose that (H) holds. If $u^{*}_{0}$ is a hyperbolic solution of (\ref{eql2}), then $u^{*}_{0}$ is isolated.
\end{corollary}

\begin{proposition}
\label{finite}
Suppose that (H) holds. If all points in $\mathcal{E}_{0}$ are isolated, then there is only a finite number of them. Moreover, if $0\notin \sigma(\Lambda-D\tilde{F}_{0}(u^{*}_{0}))$ for each $u^{*}_{0}\in \mathcal{E}_{0}$, then $\mathcal{E}_{0}$ is a finite set.
\end{proposition}
\noindent {\bf Proof. }
We suppose that the number of elements in $\mathcal{E}_{0}$ is infinite, hence there exists a sequence $\{u_{n}\}_{n\in \mathbb{N}}$ in $\mathcal{E}_{0}$. From Lemma \ref{compacto}, $\mathcal{E}_{0}$ is compact, thus there exist a subsequence $\{u_{n_{k}}\}_{k\in \mathbb{N}}$ of $\{u_{n}\}_{n\in \mathbb{N}}$ and $u^{*}\in \mathcal{E}_{0}$ such that
$$
u_{n_{k}} \to u^{*}\quad \mbox{in}\quad H^{1}(\Omega), \quad \mbox{as $k\to \infty$}.
$$ 

Thus, for all $\delta>0$, there exists $k_{0}\in \mathbb{N}$ such that $u_{n_{k}}\in B_{\delta}(u^{*})$, for all $k>k_{0}$, which is a contradiction with the fact that each fixed point in $\mathcal{E}_{0}$ is isolated and $u^{*}\in \mathcal{E}_{0}$.

Now, if $0\notin\sigma(\Lambda-D\tilde{F}_{0}(u^{*}_{0}))$ for each $u^{*}_{0}\in \mathcal{E}_{0}$, then, by Theorem \ref{theoremisolated}, $u^{*}_{0}$ is isolated. Thus, $\mathcal{E}_{0}$ is a finite set.  \quad $\blacksquare$

\vspace{0.2cm}

To prove the lower semicontinuity of the family of solutions $\left\{\mathcal{E}_{\varepsilon}\right\}_{\varepsilon\in[0,\varepsilon_{0}]}$ at $\varepsilon=0$, we will need of the following lemmas:  

\begin{lemma}
\label{lemma1lower}
Suppose that (H) holds and let $u^{*} \in H^{1}(\Omega)$. Then, for each $\varepsilon\in [0,\varepsilon_{0}]$ fixed, the operator $\Lambda^{-1}D\tilde{F}_{\varepsilon}(u^{*}): H^{1}(\Omega)\to H^{1}(\Omega)$ is compact. For any bounded family $\left\{w_{\varepsilon}\right\}_{\varepsilon\in(0,\varepsilon_{0}]}$ in $H^{1}(\Omega)$, the family $\{\Lambda^{-1}D\tilde{F}_{\varepsilon}(u^{*})w_{\varepsilon}\}_{\varepsilon\in(0,\varepsilon_{0}]}$ is relatively compact in $H^{1}(\Omega)$. Moreover, if  $w_{\varepsilon}\to w$ in $H^{1}(\Omega)$, as $\varepsilon\to0$, then 
$$
\Lambda^{-1}D\tilde{F}_{\varepsilon}(u^{*})w_{\varepsilon}\to \Lambda^{-1}D\tilde{F}_{0}(u^{*})w \quad \mbox{in}\quad H^{1}(\Omega),  \quad \mbox{as $\varepsilon\to0$}.
$$
\end{lemma}
\noindent {\bf Proof. }
For each $\varepsilon\in[0,\varepsilon_{0}]$ fixed, the compactness of linear operator $\Lambda^{-1}D\tilde{F}_{\varepsilon}(u^{*}): H^{1}(\Omega)\to H^{1}(\Omega)$ follows from item {\it 3.)} of Lemma \ref{resultsconvnonlinearity} and of compactness of linear operator $\Lambda^{-1}:H^{-s}(\Omega)\to H^{1}(\Omega)$.

Let $\{w_{\varepsilon}\}_{\varepsilon\in(0,\varepsilon_{0}]}$ be a bounded family in $H^{1}(\Omega)$. Since
$$
\|D\tilde{F}_{\varepsilon}(u^{*})w_{\varepsilon}\|_{H^{-s}(\Omega)}\leqslant \|D\tilde{F}_{\varepsilon}(u^{*})\|_{\mathcal{L}(H^{1}(\Omega),H^{-s}(\Omega))}\left\|w_{\varepsilon}\right\|_{H^{1}(\Omega)}, \quad \mbox{$\forall \varepsilon\in (0,\varepsilon_{0}]$},
$$
and from item {\it 3.)} of Lemma \ref{resultsconvnonlinearity}, $\{D\tilde{F}_{\varepsilon}(u^{*})\}_{\varepsilon\in (0,\varepsilon_{0}]}$ is a bounded family in $\mathcal{L}(H^{1}(\Omega),H^{-s}(\Omega))$, uniformly in $\varepsilon$, then $\{D\tilde{F}_{\varepsilon}(u^{*})w_{\varepsilon}\}_{\varepsilon\in (0,\varepsilon_{0}]}$ is a bounded family in $H^{-s}(\Omega)$. By compactness of $\Lambda^{-1}:H^{-s}(\Omega)\to H^{1}(\Omega)$, we have that $\{\Lambda^{-1}D\tilde{F}_{\varepsilon}(u^{*})w_{\varepsilon}\}_{\varepsilon\in (0,\varepsilon_{0}]}$ has a convergent subsequence in $H^{1}(\Omega)$. Therefore, the family $\{\Lambda^{-1}D\tilde{F}_{\varepsilon}(u^{*})w_{\varepsilon}\}_{\varepsilon\in (0,\varepsilon_{0}]}$ is relatively compact.

Now, let us take $w_{\varepsilon}\to w$ in $H^{1}(\Omega)$, as $\varepsilon\to0$. Thus, from item {\it 6.)} of Lemma \ref{resultsconvnonlinearity},
$$
D\tilde{F}_{\varepsilon}(u^{*})w_{\varepsilon}\to D\tilde{F}_{0}(u^{*})w \quad \mbox{in}\quad H^{-s}(\Omega), \quad \mbox{as $\varepsilon\to0$}.
$$  
By continuity of operator $\Lambda^{-1}:H^{-s}(\Omega)\to H^{1}(\Omega)$, we get 
$$
\Lambda^{-1}D\tilde{F}_{\varepsilon}(u^{*})w_{\varepsilon}\to \Lambda^{-1}D\tilde{F}_{0}(u^{*})w \quad \mbox{in}\quad H^{1}(\Omega), \quad \mbox{as $\varepsilon\to0$}. \quad \blacksquare
$$

\begin{lemma}\label{lemma2lower}
Suppose that (H) holds and let $u^{*}\in H^{1}(\Omega)$ such that $0\notin\sigma(\Lambda-D\tilde{F}_{0}(u^{*}))$. Then, there exist $\varepsilon_{0}>0$ and $C>0$ independent of $\varepsilon$ such that $0\notin \sigma(\Lambda-D\tilde{F}_{\varepsilon}(u^{*}))$ and 
\begin{eqnarray}
\label{resolvente1}
\|(\Lambda-D\tilde{F}_{\varepsilon}(u^{*}))^{-1}\|_{\mathcal{L}(H^{-s}(\Omega),H^{1}(\Omega))}\leqslant C,\quad \mbox {$\forall \varepsilon\in [0,\varepsilon_{0}]$}.
\end{eqnarray}
Furthermore, for each $\varepsilon\in[0,\varepsilon_{0}]$ fixed, the operator $(\Lambda-D\tilde{F}_{\varepsilon}(u^{*}))^{-1}: H^{-s}(\Omega)\to H^{1}(\Omega)$ is compact. For any bounded family $\left\{w_{\varepsilon}\right\}_{\varepsilon\in(0,\varepsilon_{0}]}$ in $H^{-s}(\Omega)$, the family $\{(\Lambda-D\tilde{F}_{\varepsilon}(u^{*}))^{-1}w_{\varepsilon}\}_{\varepsilon\in(0,\varepsilon_{0}]}$ is relatively compact in $H^{1}(\Omega)$. Moreover, if  $w_{\varepsilon}\to w$ in $H^{-s}(\Omega)$, as $\varepsilon\to0$, then 
$$
(\Lambda-D\tilde{F}_{\varepsilon}(u^{*}))^{-1}w_{\varepsilon}\to (\Lambda-D\tilde{F}_{0}(u^{*}))^{-1}w\quad \mbox{in}\quad H^{1}(\Omega), \quad \mbox{as $\varepsilon\to0$}.
$$
\end{lemma}
\noindent {\bf Proof. }
Initially, for each $\varepsilon\in[0,\varepsilon_{0}]$, we note that
$$
(\Lambda-D\tilde{F}_{\varepsilon}(u^{*}))^{-1}=[\Lambda(I-\Lambda^{-1}D\tilde{F}_{\varepsilon}(u^{*}))]^{-1}=(I-\Lambda^{-1}D\tilde{F}_{\varepsilon}(u^{*}))^{-1}\Lambda^{-1}. 
$$
Then, prove that $0\notin \sigma(\Lambda-D\tilde{F}_{\varepsilon}(u^{*}))$ it is equivalent to prove that $1\in \rho(\Lambda^{-1}D\tilde{F}_{\varepsilon}(u^{*}))$. Moreover, to prove that there exist $\varepsilon_{0}>0$ and $C>0$ independent of $\varepsilon$ such that (\ref{resolvente1}) holds, it is enough to prove that there exist $\varepsilon_{0}>0$ and $M>0$ independent of $\varepsilon$ such that
\begin{eqnarray}
\label{resolvente2}
\|(I-\Lambda^{-1}D\tilde{F}_{\varepsilon}(u^{*}))^{-1}\|_{\mathcal{L}(H^{1}(\Omega))}\leqslant M,\quad  \mbox{$\forall \varepsilon\in[0,\varepsilon_{0}]$}.
\end{eqnarray}
In fact, we note that 
$$
\begin{array}{lll}
\displaystyle\|(\Lambda-D\tilde{F}_{\varepsilon}(u^{*}))^{-1}\|_{\mathcal{L}(H^{-s}(\Omega),H^{1}(\Omega))}\\

\leqslant  \displaystyle \|(I-\Lambda^{-1}D\tilde{F}_{\varepsilon}(u^{*}))^{-1}\|_{\mathcal{L}(H^{1}(\Omega))}\|\Lambda^{-1}\|_{\mathcal{L}(H^{-s}(\Omega),H^{1}(\Omega))}\\

\leqslant \displaystyle M\|\Lambda^{-1}\|_{\mathcal{L}(H^{-s}(\Omega),H^{1}(\Omega))}=C,\quad \mbox{$\forall \varepsilon\in[0,\varepsilon_{0}]$,}
\end{array}
$$
where $C>0$ does not depend of $\varepsilon$.

Then we will show (\ref{resolvente2}). Initially, from hypothesis $0\notin\sigma(\Lambda-D\tilde{F}_{0}(u^{*}))$, then $1\in \rho(\Lambda^{-1}D\tilde{F}_{0}(u^{*}))$. Hence, there exists the inverse $(I-\Lambda^{-1}D\tilde{F}_{0}(u^{*}))^{-1}:H^{1}(\Omega)\to H^{1}(\Omega)$ and, in particular, the kernel $\mathcal{N}(I-\Lambda^{-1}D\tilde{F}_{0}(u^{*}))=\{0\}$.

Now, let $B_{\varepsilon}=\Lambda^{-1}D\tilde{F}_{\varepsilon}(u^{*})$, for all $\varepsilon\in[0,\varepsilon_{0}]$. From Lemma \ref{lemma1lower} we have that, for each $\varepsilon\in [0,\varepsilon_{0}]$ fixed, the operator $B_{\varepsilon}: H^{1}(\Omega)\to H^{1}(\Omega)$ is compact. Using the compactness of $B_{\varepsilon}$, we can show that estimate (\ref{resolvente2}) is equivalent to say
\begin{eqnarray}
\label{resolvente3}
\|(I-B_{\varepsilon})u_{\varepsilon}\|_{H^{1}(\Omega)}\geqslant \frac{1}{M},\quad  \mbox{$\forall \varepsilon\in[0,\varepsilon_{0}]$}\quad \mbox{and}\quad \mbox{$\forall  u_{\varepsilon}\in H^{1}(\Omega)$\quad with\quad $\left\|u_{\varepsilon}\right\|_{H^{1}(\Omega)}=1$}.
\end{eqnarray}

In fact, suppose that (\ref{resolvente2}) holds, then there exists the inverse $(I-B_{\varepsilon})^{-1}:H^{1}(\Omega)\to H^{1}(\Omega)$ and it is continuous. Moreover,
$$
\|(I-B_{\varepsilon})^{-1}v_{\varepsilon}\|_{H^{1}(\Omega)}\leqslant M\left\|v_{\varepsilon}\right\|_{H^{1}(\Omega)},\quad \mbox{$\forall \varepsilon\in[0,\varepsilon_{0}]$}\quad \mbox{and}\quad \mbox{$\forall v_{\varepsilon}\in H^{1}(\Omega)$}.
$$ 

Let $u_{\varepsilon}\in H^{1}(\Omega)$ such that $\left\|u_{\varepsilon}\right\|_{H^{1}(\Omega)}=1$ and taking $v_{\varepsilon}=(I-B_{\varepsilon})u_{\varepsilon}$, we have
$$
\|(I-B_{\varepsilon})^{-1}(I-B_{\varepsilon})u_{\varepsilon}\|_{H^{1}(\Omega)}\leqslant M \|(I-B_{\varepsilon})u_{\varepsilon}\|_{H^{1}(\Omega)}
$$
$$
\Rightarrow 1=\left\|u_{\varepsilon}\right\|_{H^{1}(\Omega)}\leqslant M\|(I-B_{\varepsilon})u_{\varepsilon}\|_{H^{1}(\Omega)}\Rightarrow \|(I-B_{\varepsilon})u_{\varepsilon}\|_{H^{1}(\Omega)}\geqslant \frac{1}{M}.
$$

Therefore, (\ref{resolvente3}) holds. Reversely, suppose that (\ref{resolvente3}) holds. We want to prove that there exists the inverse $(I-B_{\varepsilon})^{-1}:H^{1}(\Omega)\to H^{1}(\Omega)$, it is continuous and satisfies (\ref{resolvente2}). For this, we will prove the following estimative
\begin{eqnarray}
\label{resolvente4}
\|(I-B_{\varepsilon})u_{\varepsilon}\|_{H^{1}(\Omega)}\geqslant \frac{1}{M}\left\|u_{\varepsilon}\right\|_{H^{1}(\Omega)},\quad  \mbox{$\forall \varepsilon\in[0,\varepsilon_{0}]$}\quad \mbox{and}\quad \mbox{$\forall u_{\varepsilon}\in H^{1}(\Omega)$}.
\end{eqnarray}

We note that (\ref{resolvente4}) is immediate for $u_{\varepsilon}=0$. Let $u_{\varepsilon}\in H^{1}(\Omega)$, $u_{\varepsilon}\neq0$, and we take  $v_{\varepsilon}=\displaystyle\frac{u_{\varepsilon}}{\left\|u_{\varepsilon}\right\|_{H^{1}(\Omega)}}$. Thus, $\left\|v_{\varepsilon}\right\|_{H^{1}(\Omega)}=1$ and using (\ref{resolvente3}), we get  
$$
\|(I-B_{\varepsilon})v_{\varepsilon}\|_{H^{1}(\Omega)}\geqslant \frac{1}{M}\Rightarrow \left\|(I-B_{\varepsilon})\frac{u_{\varepsilon}}{\left\|u_{\varepsilon}\right\|_{H^{1}(\Omega)}}\right\|_{H^{1}(\Omega)}\geqslant \frac{1}{M} 
$$
$$
\Rightarrow \frac{1}{\left\|u_{\varepsilon}\right\|_{H^{1}(\Omega)}} \|(I-B_{\varepsilon})u_{\varepsilon}\|_{H^{1}(\Omega)}\geqslant \frac{1}{M}\Rightarrow \|(I-B_{\varepsilon})u_{\varepsilon}\|_{H^{1}(\Omega)}\geqslant \frac{1}{M} \left\|u_{\varepsilon}\right\|_{H^{1}(\Omega)}.
$$

Now, let $u_{\varepsilon}\in H^{1}(\Omega)$ such that $(I-B_{\varepsilon})u_{\varepsilon}=0$. From (\ref{resolvente4}) follows $u_{\varepsilon}=0$. Thus, for each $\varepsilon\in[0,\varepsilon_{0}]$, $\mathcal{N}(I-B_{\varepsilon})=\{0\}$ and the operator $I-B_{\varepsilon}$ is injective. So there exists the inverse $(I-B_{\varepsilon})^{-1}:\mathcal{R}(I-B_{\varepsilon})\to H^{1}(\Omega)$, where $\mathcal{R}(I-B_{\varepsilon})$ denotes the image of the operator $I-B_{\varepsilon}$. 

Since $B_{\varepsilon}$ is compact, for all $\varepsilon\in[0,\varepsilon_{0}]$, then by Fredholm Alternative Theorem, we have 
$$
\mathcal{N}(I-B_{\varepsilon})=\{0\}\Leftrightarrow \mathcal{R}(I-B_{\varepsilon})=H^{1}(\Omega).
$$ 
Hence, $\mathcal{R}(I-B_{\varepsilon})=H^{1}(\Omega)$ and $I-B_{\varepsilon}$ is bijective, thus there exists the inverse $(I-B_{\varepsilon})^{-1}:H^{1}(\Omega)\to H^{1}(\Omega)$.

Now, taking $v_{\varepsilon}\in H^{1}(\Omega)$ we have that there exists $u_{\varepsilon}\in H^{1}(\Omega)$ such that $(I-B_{\varepsilon})u_{\varepsilon}=v_{\varepsilon}$ and $u_{\varepsilon}=(I-B_{\varepsilon})^{-1}v_{\varepsilon}$. From (\ref{resolvente4}) we have
$$
\|(I-B_{\varepsilon})^{-1}v_{\varepsilon}\|_{H^{1}(\Omega)}=\left\|u_{\varepsilon}\right\|_{H^{1}(\Omega)}\leqslant M\|(I-B_{\varepsilon})u_{\varepsilon}\|_{H^{1}(\Omega)}=M\left\|v_{\varepsilon}\right\|_{H^{1}(\Omega)}
$$ 
$$
\Rightarrow \|(I-B_{\varepsilon})^{-1}\|_{\mathcal{L}(H^{1}(\Omega))} \leqslant M,\quad  \mbox{$\forall \varepsilon\in[0,\varepsilon_{0}]$}. 
$$
Therefore, (\ref{resolvente2}) holds.

Since (\ref{resolvente2}) and (\ref{resolvente3}) are equivalents, then it is enough to show (\ref{resolvente3}). Suppose that (\ref{resolvente3}) is not true, that is, there exist a sequence $\left\{u_{n}\right\}_{n\in \mathbb{N}}$ in $H^{1}(\Omega)$, with $\left\|u_{n}\right\|_{H^{1}(\Omega)}=1$ and $\varepsilon_{n}\to 0$, as $n\to\infty$, such that
$$
\|(I-B_{\varepsilon_{n}})u_{n}\|_{H^{1}(\Omega)}\to 0, \quad \mbox{as $n\to\infty$}.
$$

From Lemma \ref{lemma1lower} we get that $\left\{B_{\varepsilon_{n}}u_{n}\right\}_{n\in \mathbb{N}}$ is relatively compact. Thus, $\left\{B_{\varepsilon_{n}}u_{n}\right\}_{n\in \mathbb{N}}$ has a convergent subsequence, which we again denote by $\left\{B_{\varepsilon_{n}}u_{n}\right\}_{n\in \mathbb{N}}$, with limit $u\in H^{1}(\Omega)$, that is,
$$
B_{\varepsilon_{n}}u_{n}\to u \quad \mbox{in}\quad H^{1}(\Omega), \quad \mbox{as $n\to \infty$}.
$$ 

Since $u_{n}-B_{\varepsilon_{n}}u_{n}\to 0$ in $H^{1}(\Omega)$, as $n\to\infty$, then $u_{n}\to u$ in $H^{1}(\Omega)$, as $n\to\infty$. Hence, $\left\|u\right\|_{H^{1}(\Omega)}=1$. Moreover, since $u_{n}\to u$ in $H^{1}(\Omega)$, as $n\to\infty$, then using the Lemma \ref{lemma1lower}, we have $B_{\varepsilon_{n}}u_{n}\to B_{0}u$ in $H^{1}(\Omega)$, as $n\to\infty$. Thus,
$$
u_{n}-B_{\varepsilon_{n}}u_{n}\to u-B_{0}u\quad \mbox{in}\quad H^{1}(\Omega), \quad \mbox{as $n\to\infty$}.
$$

By the uniqueness of the limit, $u-B_{0}u=0$. This implies that $(I-B_{0})u=0$, with $u\neq0$, contradicting our hypothesis. Therefore, (\ref{resolvente3}) holds.

With this we conclude that there exist $\varepsilon_{0}>0$ and $C>0$ independent of $\varepsilon$ such that (\ref{resolvente1}) holds.

Now, for each $\varepsilon\in[0,\varepsilon_{0}]$, the operator $(\Lambda-D\tilde{F}_{\varepsilon}(u^{*}))^{-1}$ is compact and the prove of this compactness follows similarly to account below.  

Let $\{w_{\varepsilon}\}_{\varepsilon\in(0,\varepsilon_{0}]}$ be a bounded family in $H^{-s}(\Omega)$. For each $\varepsilon\in(0,\varepsilon_{0}]$, let  $v_{\varepsilon}=(\Lambda-D\tilde{F}_{\varepsilon}(u^{*}))^{-1}w_{\varepsilon}.$ From (\ref{resolvente1}) we have
$$
\begin{array}{lll}
\displaystyle \left\|v_{\varepsilon}\right\|_{H^{1}(\Omega)} &\leqslant& \|(\Lambda-D\tilde{F}_{\varepsilon}(u^{*}))^{-1}w_{\varepsilon}\|_{H^{1}(\Omega)} \\

\displaystyle &\leqslant& \|(\Lambda-D\tilde{F}_{\varepsilon}(u^{*}))^{-1}\|_{\mathcal{L}(H^{-s}(\Omega),H^{1}(\Omega))}\left\|w_{\varepsilon}\right\|_{H^{-s}(\Omega)}\leqslant C\left\|w_{\varepsilon}\right\|_{H^{-s}(\Omega)}.
\end{array}
$$
Hence, $\{v_{\varepsilon}\}_{\varepsilon\in(0,\varepsilon_{0}]}$ is a bounded family in $H^{1}(\Omega)$. Moreover, 
$$
v_{\varepsilon}=(\Lambda-D\tilde{F}_{\varepsilon}(u^{*}))^{-1}w_{\varepsilon}=(I-\Lambda^{-1}D\tilde{F}_{\varepsilon}(u^{*}))^{-1}\Lambda^{-1}w_{\varepsilon}
$$
$$
\Leftrightarrow (I-\Lambda^{-1}D\tilde{F}_{\varepsilon}(u^{*}))v_{\varepsilon}=\Lambda^{-1}w_{\varepsilon}\Leftrightarrow v_{\varepsilon}=\Lambda^{-1}D\tilde{F}_{\varepsilon}(u^{*})v_{\varepsilon}+\Lambda^{-1}w_{\varepsilon}.
$$

By compactness of $\Lambda^{-1}:H^{-s}(\Omega)\to H^{1}(\Omega)$, we get that $\{\Lambda^{-1}w_{\varepsilon}\}_{\varepsilon\in(0,\varepsilon_{0}]}$ has a convergent subsequence in $H^{1}(\Omega)$. Moreover, using the Lemma \ref{lemma1lower}, we have that $\{\Lambda^{-1}D\tilde{F}_{\varepsilon}(u^{*})v_{\varepsilon}\}_{\varepsilon\in(0,\varepsilon_{0}]}$ is relatively compact in $H^{1}(\Omega)$, then $\{\Lambda^{-1}D\tilde{F}_{\varepsilon}(u^{*})v_{\varepsilon}\}_{\varepsilon\in(0,\varepsilon_{0}]}$ has a convergent subsequence in $H^{1}(\Omega)$. Therefore, $\{v_{\varepsilon}\}_{\varepsilon\in(0,\varepsilon_{0}]}$ has a convergent subsequence in $H^{1}(\Omega)$, that is, the family  $\{(\Lambda-D\tilde{F}_{\varepsilon}(u^{*}))^{-1}w_{\varepsilon}\}_{\varepsilon\in(0,\varepsilon_{0}]}$ has a convergent subsequence in $H^{1}(\Omega)$, thus it is relatively compact in $H^{1}(\Omega)$.

Now, we take $w_{\varepsilon}\to w$ in $H^{-s}(\Omega)$, as $\varepsilon\to0$. By continuity of operator $\Lambda^{-1}:H^{-s}(\Omega)\to H^{1}(\Omega)$, we have
$$
\Lambda^{-1}w_{\varepsilon}\to \Lambda^{-1}w \quad \mbox{in}\quad H^{1}(\Omega), \quad \mbox{as $\varepsilon\to0$}.
$$
Moreover, $\{w_{\varepsilon}\}_{\varepsilon\in(0,\varepsilon_{0}]}$ is bounded in $H^{-s}(\Omega)$, for some $\varepsilon_{0}>0$ sufficiently small, and we have that from the above that $\{v_{\varepsilon}\}_{\varepsilon\in(0,\varepsilon_{0}]}$, with $\varepsilon_{0}>0$ sufficiently small, has a convergent subsequence, which we again denote by $\{v_{\varepsilon}\}_{\varepsilon\in(0,\varepsilon_{0}]}$, with limit $v\in H^{1}(\Omega)$, that is, 
$$
v_{\varepsilon}\to v\quad \mbox{in}\quad H^{1}(\Omega),\quad \mbox{as $\varepsilon\to0$}.
$$ 

From Lemma \ref{lemma1lower} we get
$$
\Lambda^{-1}D\tilde{F}_{\varepsilon}(u^{*})v_{\varepsilon}\to \Lambda^{-1}D\tilde{F}_{0}(u^{*})v\quad \mbox{in}\quad H^{1}(\Omega), \quad \mbox{as $\varepsilon\to0$}.
$$ 
Thus, $v$ satisfies $v=\Lambda^{-1}D\tilde{F}_{0}(u^{*})v+\Lambda^{-1}w$, and so $v=(\Lambda-D\tilde{F}_{0}(u^{*}))^{-1}w$. Therefore,
$$
(\Lambda-D\tilde{F}_{\varepsilon}(u^{*}))^{-1}w_{\varepsilon}\to (\Lambda-D\tilde{F}_{0}(u^{*}))^{-1}w \quad \mbox{in}\quad H^{1}(\Omega),\quad \mbox{as $\varepsilon\to0$}.
$$
The limit above is independent of the subsequence, thus whole family $\{(\Lambda-D\tilde{F}_{\varepsilon}(u^{*}))^{-1}w_{\varepsilon}\}_{\varepsilon\in(0,\varepsilon_{0}]}$ converges to $(\Lambda-D\tilde{F}_{0}(u^{*}))^{-1}w$ in $H^{1}(\Omega)$, as $\varepsilon\to 0$. \quad $\blacksquare$

\begin{theorem}
\label{theolowerequilibria}
Suppose that (H) holds and that $u^{*}_{0}$ is a solution of (\ref{eql2}) which satisfies $0\notin\sigma(\Lambda-D\tilde{F}_{0}(u^{*}_{0}))$. Then, there exist $\varepsilon_{0}>0$ and $\delta>0$ such that, for each $\varepsilon \in (0, \varepsilon_{0}]$, the equation (\ref{eql1}) has exactly one solution, $u^{*}_{\varepsilon}$, in 
$$
\{v_{\varepsilon}\in H^{1}(\Omega):  \left\|v_{\varepsilon}-u^{*}_{0}\right\|_{H^{1}(\Omega)}\leqslant \delta\}.
$$
Furthermore,
$$
u^{*}_{\varepsilon}\to u^{*}_{0}\quad \mbox{in}\quad H^{1}(\Omega), \quad \mbox{as $\varepsilon\to0$}.
$$
In particular, the family of solutions $\left\{\mathcal{E}_{\varepsilon}\right\}_{\varepsilon\in[0,\varepsilon_{0}]}$ of (\ref{eql1}) and (\ref{eql2}) is lower semicontinuous at $\varepsilon=0$.
\end{theorem}
\noindent {\bf Proof. }
Initially, using the Lemma \ref{lemma2lower}, we have that there exist $\varepsilon_{0}>0$ and $C>0$, independent of $\varepsilon$, such that $0\notin\sigma(\Lambda-D\tilde{F}_{\varepsilon}(u^{*}_{0}))$ and
\begin{eqnarray}
\label{deslem3eq}
\|(\Lambda-D\tilde{F}_{\epsilon}(u^{*}_{0}))^{-1}\|_{\mathcal{L}(H^{-s}(\Omega),H^{1}(\Omega))}\leqslant
C,\quad  \mbox{$\forall \varepsilon\in(0,\varepsilon_{0}]$}.
\end{eqnarray}

By item {\it 1.)} of Lemma \ref{resultsconvnonlinearity}, there exists $\tilde{\delta}=\tilde{\delta}(C)>0$ independent of $\varepsilon$ such that  
\begin{eqnarray}
\label{deslem7eq}
C\|\tilde{F}_{\varepsilon}(u_{\varepsilon})-\tilde{F}_{\varepsilon}(v_{\varepsilon})-D\tilde{F}_{\varepsilon}(u^{*}_{0})(u_{\varepsilon}-v_{\varepsilon})\|_{H^{-s}(\Omega)}\leqslant \frac{1}{2}\left\|u_{\varepsilon}-v_{\varepsilon}\right\|_{H^{1}(\Omega)}, \quad \mbox{$\forall  \varepsilon\in(0,\varepsilon_{0}]$},
\end{eqnarray}
for $\left\|u_{\varepsilon}-v_{\varepsilon}\right\|\leqslant \tilde{\delta}$.

We note that $u_{\varepsilon}$, $\varepsilon \in (0, \varepsilon_{0}]$, is a solution of (\ref{eql1}) if and only if $u_{\varepsilon}$ is a fixed point of the map
$$
\begin{array}{ll}
\Phi_{\varepsilon}: H^{1}(\Omega)\to H^{1}(\Omega)\\
\qquad \qquad u_{\varepsilon}\mapsto \Phi_{\varepsilon}(u_{\varepsilon})=(\Lambda-D\tilde{F}_{\varepsilon}(u^{*}_{0}))^{-1}(\tilde{F}_{\varepsilon}(u_{\varepsilon})-D\tilde{F}_{\varepsilon}(u^{*}_{0})u_{\varepsilon}).
\end{array}
$$

Initially, we affirm that 
\begin{eqnarray}
\label{desconv1eq}
\Phi_{\varepsilon}(u^{*}_{0})\to u^{*}_{0}\quad \mbox{in}\quad H^{1}(\Omega), \quad \mbox{as $\varepsilon\to0$}.
\end{eqnarray}
In fact, using (\ref{deslem3eq}), for $\varepsilon \in (0, \varepsilon_{0}]$, we have
$$
\|\Phi_{\varepsilon}(u^{*}_{0})-u^{*}_{0}\|_{H^{1}(\Omega)}
$$
$$
\begin{array}{rcl}
&\leqslant& \displaystyle \| (\Lambda-D\tilde{F}_{\varepsilon}(u^{*}_{0}))^{-1}[(\tilde{F}_{\varepsilon}(u^{*}_{0})-D\tilde{F}_{\varepsilon}(u^{*}_{0})u^{*}_{0})-(\tilde{F}_{0}(u^{*}_{0})-D\tilde{F}_{0}(u^{*}_{0})u^{*}_{0})]\|_{H^{1}(\Omega)}\\

&+& \displaystyle \| [ (\Lambda-D\tilde{F}_{\varepsilon}(u^{*}_{0}))^{-1}- (\Lambda-D\tilde{F}_{0}(u^{*}_{0}))^{-1}] (\tilde{F}_{0}(u^{*}_{0})-D\tilde{F}_{0}(u^{*}_{0})u^{*}_{0}) \|_{H^{1}(\Omega)}\\

&\leqslant& \displaystyle  C( \|\tilde{F}_{\varepsilon}(u^{*}_{0})-\tilde{F}_{0}(u^{*}_{0})\|_{H^{-s}(\Omega)}+\|D\tilde{F}_{\varepsilon}(u^{*}_{0})u^{*}_{0}-D\tilde{F}_{0}(u^{*}_{0})u^{*}_{0}\|_{H^{-s}(\Omega)} )\\

&+&\displaystyle \| [ (\Lambda-D\tilde{F}_{\varepsilon}(u^{*}_{0}))^{-1}- (\Lambda-D\tilde{F}_{0}(u^{*}_{0}))^{-1}] (\tilde{F}_{0}(u^{*}_{0})-D\tilde{F}_{0}(u^{*}_{0})u^{*}_{0}) \|_{H^{1}(\Omega)}\to0,
\end{array}
$$
as $\varepsilon\to0$. This follows from \cite[Lemma 2.1]{aragaobezerra1}, item {\it 6.)} of Lemma \ref{resultsconvnonlinearity} and Lemma \ref{lemma2lower}.

Next, we show that, for $\varepsilon \in (0, \varepsilon_{0}]$, for some $\varepsilon_{0}>0$ sufficiently small, $\Phi_{\varepsilon}$ is a contraction map from the closed ball $\bar{B}_{\delta}(u^{*}_{0})=\{v_{\varepsilon}\in H^{1}(\Omega): \left\|v_{\varepsilon}-u^{*}_{0}\right\|_{H^{1}(\Omega)}\leqslant \delta\}$ into itself, where $\delta=\frac{\tilde{\delta}}{2}$. First, we show that $\Phi_{\varepsilon}$ is a contraction on the $\bar{B}_{\delta}(u^{*}_{0})$ (uniformly in $\varepsilon$). Let $u_{\varepsilon}$, $v_{\varepsilon}\in \bar{B}_{\delta}(u^{*}_{0})$ and using (\ref{deslem3eq}) and (\ref{deslem7eq}), for $\varepsilon \in (0, \varepsilon_{0}]$, we have
$$
\begin{array}{lll}
\displaystyle\left\|\Phi_{\varepsilon}(u_{\varepsilon})-\Phi_{\varepsilon}(v_{\varepsilon}) \right\|_{H^{1}(\Omega)}\\

= \displaystyle \|(\Lambda-D\tilde{F}_{\varepsilon}(u^{*}_{0}))^{-1}[\tilde{F}_{\varepsilon}(u_{\varepsilon})-\tilde{F}_{\varepsilon}(v_{\varepsilon})-D\tilde{F}_{\varepsilon}(u^{*}_{0})(u_{\varepsilon}-v_{\varepsilon})]\|_{H^{1}(\Omega)}\\

\leqslant  \displaystyle C  \|\tilde{F}_{\varepsilon}(u_{\varepsilon})-\tilde{F}_{\varepsilon}(v_{\varepsilon})-D\tilde{F}_{\varepsilon}(u^{*}_{0})(u_{\varepsilon}-v_{\varepsilon})\|_{H^{-s}(\Omega)}\\

\leqslant  \displaystyle \frac{1}{2}\left\|u_{\varepsilon}-v_{\varepsilon}\right\|_{H^{1}(\Omega)},\quad \mbox{for $\varepsilon \in (0,\varepsilon_{0}]$.}
\end{array} 
$$
To show that $\Phi_{\varepsilon}$ maps $\bar{B}_{\delta}(u^{*}_{0})$ into itself, we observe that if $u_{\varepsilon}\in \bar{B}_{\delta}(u^{*}_{0})$, then 
$$
\begin{array}{rcl}
\displaystyle\left\|\Phi_{\varepsilon}(u_{\varepsilon})-u^{*}_{0}\right\|_{H^{1}(\Omega)}&\leqslant& \displaystyle \left\|\Phi_{\varepsilon}(u_{\varepsilon})-\Phi_{\varepsilon}(u^{*}_{0}) \right\|_{H^{1}(\Omega)}+\left\|\Phi_{\varepsilon}(u^{*}_{0})-u^{*}_{0}\right\|_{H^{1}(\Omega)}\\

&\leqslant& \displaystyle \frac{\delta}{2}+\left\|\Phi_{\varepsilon}(u^{*}_{0})-u^{*}_{0}\right\|_{H^{1}(\Omega)}, \quad \mbox{for $\varepsilon\in (0,\varepsilon_{0}]$.} 
\end{array}
$$
By convergence in (\ref{desconv1eq}), we have that there exists $\varepsilon_{0}>0$ such that
$$
\left\|\Phi_{\varepsilon}(u_{\varepsilon})-u^{*}_{0}\right\|_{H^{1}(\Omega)}\leqslant \frac{\delta}{2}+\frac{\delta}{2}=\delta, \quad \mbox{for $\varepsilon\in(0,\varepsilon_{0}]$}.
$$
Hence, $\Phi_{\varepsilon}:\bar{B}_{\delta}(u^{*}_{0})\to \bar{B}_{\delta}(u^{*}_{0})$ is a contraction for all $\varepsilon \in (0, \varepsilon_{0}]$. By Contraction Theorem follows that, for each $\varepsilon \in (0, \varepsilon_{0}]$, $\Phi_{\varepsilon}$ has an unique fixed point, $u^{*}_{\varepsilon}$, in the $\bar{B}_{\delta}(u^{*}_{0})$. 

To show that $u^{*}_{\varepsilon}\to u^{*}_{0}$ in $H^{1}(\Omega)$, as $\varepsilon\to0$, we proceed in the following manner: since $\Phi_{\varepsilon}$ is a contraction map from $\bar{B}_{\delta}(u^{*}_{0})$ into itself, then
$$
\begin{array}{lll}
\displaystyle\left\|u^{*}_{\varepsilon}-u^{*}_{0}\right\|_{H^{1}(\Omega)}&=& \displaystyle\left\|\Phi_{\varepsilon}(u^{*}_{\varepsilon})-u^{*}_{0}\right\|_{H^{1}(\Omega)}\\

&\leqslant& \displaystyle \left\|\Phi_{\varepsilon}(u^{*}_{\varepsilon})-\Phi_{\varepsilon}(u^{*}_{0})\right\|_{H^{1}(\Omega)}+\left\|\Phi_{\varepsilon}(u^{*}_{0})-u^{*}_{0}\right\|_{H^{1}(\Omega)}\\

&\leqslant& \displaystyle \frac{1}{2} \left\|u^{*}_{\varepsilon}-u^{*}_{0}\right\|_{H^{1}(\Omega)}+\left\|\Phi_{\varepsilon}(u^{*}_{0})-u^{*}_{0}\right\|_{H^{1}(\Omega)}.
\end{array}
$$
Thus, using (\ref{desconv1eq}),
$$
\left\|u^{*}_{\varepsilon}-u^{*}_{0}\right\|_{H^{1}(\Omega)}\leqslant 2\left\|\Phi_{\varepsilon}(u^{*}_{0})-u^{*}_{0}\right\|_{H^{1}(\Omega)}\to 0, \quad \mbox{as $\varepsilon\to0$}. 
$$
Hence and by compactness of $\mathcal{E}_{0}$ (Lemma \ref{compacto}), we have that the family $\left\{\mathcal{E}_{\varepsilon}\right\}_{\varepsilon\in[0,\varepsilon_{0}]}$ is lower semicontinuity at $\varepsilon=0$. \quad $\blacksquare$

\begin{corollary}
\label{{corollarylowerequilibria}}
Suppose that (H) holds and that $u^{*}_{0}$ is a hyperbolic solution of (\ref{eql2}). Then, there exist $\varepsilon_{0}>0$ and $\delta>0$ such that, for each $\varepsilon \in (0, \varepsilon_{0}]$, the equation (\ref{eql1}) has  exactly one solution, $u^{*}_{\varepsilon}$, in 
$$
\{v_{\varepsilon}\in H^{1}(\Omega): \left\|v_{\varepsilon}-u^{*}_{0}\right\|_{H^{1}(\Omega)}\leqslant \delta\}.
$$
Furthermore, $u^{*}_{\varepsilon}\to u^{*}_{0}$ in $H^{1}(\Omega)$, as $\varepsilon\to0$.
\end{corollary}

\begin{remark}
\label{remarkhyperbolic}
Now that we have obtained an unique solution $u^{*}_{\varepsilon}$ for (\ref{eql1}) in a small neighborhood of the hyperbolic solution $u^{*}_{0}$ for (\ref{eql2}), we can consider the linearization $\Lambda-D\tilde{F}_{\varepsilon}(u^{*}_{\varepsilon})$ and from the convergence of $u^{*}_{\varepsilon}$ to $u^{*}_{0}$ in $H^{1}(\Omega)$ it is easy to obtain that $(\Lambda-D\tilde{F}_{\varepsilon}(u^{*}_{\varepsilon}))^{-1}w_{\varepsilon}$ converges to $(\Lambda-D\tilde{F}_{0}(u^{*}_{0}))^{-1}w$ in $H^{1}(\Omega)$, whenever $w_{\varepsilon}\to w$ in $H^{-s}(\Omega)$, as $\varepsilon\to0$. Consequently, the hyperbolicity of $u^{*}_{0}$ implies the hyperbolicity of $u^{*}_{\varepsilon}$, for suitably small $\varepsilon$.
\end{remark}

\begin{theorem}
\label{theofiniteequilibria}
Suppose that (H) holds. If all solutions $u^{*}_{0}$ of (\ref{eql2}) satisfy $0\notin\sigma(\Lambda-D\tilde{F}_{0}(u^{*}_{0}))$, then (\ref{eql2}) has a finite number $k$ of solutions, $u^{*}_{0,1},...,u^{*}_{0,k}$, and there exists $\varepsilon_{0}>0$ such that, for each $\varepsilon \in (0, \varepsilon_{0}]$, the equation (\ref{eql1}) has exactly $k$ solutions, $u^{*}_{\varepsilon,1},...,u^{*}_{\varepsilon,k}$. Moreover, for all $i=1,...,k$,
$$
u^{*}_{\varepsilon,i}\to u^{*}_{0,i}\quad \mbox{in}\quad H^{1}(\Omega), \quad \mbox{as $\varepsilon\to0$}.
$$
\end{theorem}
\noindent {\bf Proof. }
The proof follows of Proposition \ref{finite} and Theorem \ref{theolowerequilibria}. \quad $\blacksquare$


\end{document}